%
%
%

\ifx\begin\undefined\else\global\advance\srcdepth by
1\expandafter \fi

\def\begin{}
\newcount\srcdepth
\srcdepth=1

\outer\def\bye{\global\advance\srcdepth by -1
  \ifnum\srcdepth=0
    \def\endcmd{\vfill\supereject\nopagenumbers\par\vfill\supereject\end}
  \else\def\endcmd{}\fi
  \endcmd
}



\def\initialize#1#2#3#4#5#6{
  \ifnum\srcdepth=1
  \magnification=#1
  \hsize = #2
  \vsize = #3
  \hoffset=#4
  \advance\hoffset by -\hsize
  \divide\hoffset by 2
  \advance\hoffset by -1truein
  \voffset=#5
  \advance\voffset by -\vsize
  \divide\voffset by 2
  \advance\voffset by -1truein
  \advance\voffset by #6
  \baselineskip=13pt
  \emergencystretch = 0.05\hsize
  \fi
}

\def\print{\initialize{1095}
  {5.5truein}{8.5truein}{8.5truein}{11truein}{-.0625truein}}

\newif\ifblackboardbold

\blackboardboldtrue


\font\titlefont=cmbx12 scaled\magstephalf
\font\sectionfont=cmbx12

\font\scriptit=cmti10 at 7pt
\font\scriptsl=cmsl10 at 7pt
\scriptfont\itfam=\scriptit
\scriptfont\slfam=\scriptsl


\newfam\bboldfam
\ifblackboardbold
\font\tenbbold=msbm10
\font\sevenbbold=msbm7
\font\fivebbold=msbm5
\textfont\bboldfam=\tenbbold
\scriptfont\bboldfam=\sevenbbold
\scriptscriptfont\bboldfam=\fivebbold
\def\bbold{\fam\bboldfam\tenbbold}
\else
\def\bbold{\bf}
\fi


\newfam\msamfam
\font\tenmsam=msam10
\font\sevenmsam=msam7
\font\fivemsam=msam5
\textfont\msamfam=\tenmsam
\scriptfont\msamfam=\sevenmsam
\scriptscriptfont\msamfam=\fivemsam

\newfam\msbmfam
\font\tenmsbm=msam10
\font\sevenmsbm=msam7
\font\fivemsbm=msam5
\textfont\msbmfam=\tenmsbm
\scriptfont\msbmfam=\sevenmsbm
\scriptscriptfont\msbmfam=\fivemsbm

\newcount\amsfamcount 
\newcount\classcount   
\newcount\positioncount
\newcount\codecount
\newcount\n             
\def\newsymbol#1#2#3#4#5{               
\n="#2                                  
\ifnum\n=1 \amsfamcount=\msamfam\else   
\ifnum\n=2 \amsfamcount=\msbmfam\else   
\ifnum\n=3 \amsfamcount=\eufmfam
\fi\fi\fi
\multiply\amsfamcount by "100           
\classcount="#3                 
\multiply\classcount by "1000           
\positioncount="#4#5            
\codecount=\classcount                  
\advance\codecount by \amsfamcount      
\advance\codecount by \positioncount
\mathchardef#1=\codecount}              


\font\Arm=cmr9
\font\Ai=cmmi9
\font\Asy=cmsy9
\font\Abf=cmbx9
\font\Brm=cmr7
\font\Bi=cmmi7
\font\Bsy=cmsy7
\font\Bbf=cmbx7
\font\Crm=cmr6
\font\Ci=cmmi6
\font\Csy=cmsy6
\font\Cbf=cmbx6

\ifblackboardbold
\font\Abbold=msbm10 at 9pt
\font\Bbbold=msbm7
\font\Cbbold=msbm5 at 6pt
\fi

\def\small{%
\textfont0=\Arm \scriptfont0=\Brm \scriptscriptfont0=\Crm
\textfont1=\Ai \scriptfont1=\Bi \scriptscriptfont1=\Ci
\textfont2=\Asy \scriptfont2=\Bsy \scriptscriptfont2=\Csy
\textfont\bffam=\Abf \scriptfont\bffam=\Bbf \scriptscriptfont\bffam=\Cbf
\def\rm{\fam0\Arm}\def\mit{\fam1}\def\oldstyle{\fam1\Ai}%
\def\bf{\fam\bffam\Abf}%
\ifblackboardbold
\textfont\bboldfam=\Abbold
\scriptfont\bboldfam=\Bbbold
\scriptscriptfont\bboldfam=\Cbbold
\def\bbold{\fam\bboldfam\Abbold}%
\fi
\rm
}








\newlinechar=`@
\def\forwardmsg#1#2#3{\immediate\write16{@*!*!*!* forward reference should
be: @\noexpand\forward{#1}{#2}{#3}@}}
\def\nodefmsg#1{\immediate\write16{@*!*!*!* #1 is an undefined reference@}}

\def\forwardsub#1#2{\def\newref{{#2}{#1}}}

\def\forward#1#2#3{%
\expandafter\expandafter\expandafter\forwardsub\expandafter{#3}{#2}
\expandafter\ifx\csname#1\endcsname\relax\else%
\expandafter\ifx\csname#1\endcsname\newref\else%
\forwardmsg{#1}{#2}{#3}\fi\fi%
\expandafter\let\csname#1\endcsname\newref}

\def\firstarg#1{\expandafter\argone #1}\def\argone#1#2{#1}
\def\secondarg#1{\expandafter\argtwo #1}\def\argtwo#1#2{#2}

\def\ref#1{\expandafter\ifx\csname#1\endcsname\relax
  {\nodefmsg{#1}\bf`#1'}\else
  \expandafter\firstarg\csname#1\endcsname
  ~\htmllocref{#1}{\expandafter\secondarg\csname#1\endcsname}\fi}

\def\refs#1{\expandafter\ifx\csname#1\endcsname\relax
  {\nodefmsg{#1}\bf`#1'}\else
  \expandafter\firstarg\csname #1\endcsname
  s~\htmllocref{#1}{\expandafter\secondarg\csname#1\endcsname}\fi}

\def\refn#1{\expandafter\ifx\csname#1\endcsname\relax
  {\nodefmsg{#1}\bf`#1'}\else
  \htmllocref{#1}{\expandafter\secondarg\csname #1\endcsname}\fi}



\def\widow#1{\vskip 0pt plus#1\vsize\goodbreak\vskip 0pt plus-#1\vsize}

%
\def\begincolor#1{\ifhypers\fi}%
\def\endcolor{\ifhypers\fi}%
\def\colored#1#2{%
  \begincolor{#1}#2\endcolor}
\def\red{[0.8 0.5 0]}%
\def\linkcolor{[1.0 0.1 0.1]}%
\def\sectioncolor{[0 0.4 0.8]}%
\def\proclaimcolor{[0 0.4 0.8]}%
\def\proofcolor{[0 0 0]}%
%
%
\newif\ifhypers                 
\hypersfalse
%
%
\edef\freehash{\catcode`\noexpand\#=\the\catcode`\#}%
\catcode`\#=12
\freehash
\let\freehash=\relax
\ifhypers\fi
\def\puthtml#1{\ifhypers\fi}
\def\htmlanchor#1#2{\puthtml{<a name="#1">}#2\puthtml{</a>}}

\def\@pdfm@mark#1{}
\def\setlink#1{\colored{\linkcolor}{#1}}%
%
%
%
%
\def\htmllocref#1#2{\ifhypers\leavevmode\fi\setlink{#2}\ifhypers\fi\relax}%
%
%
%
%
\def\Acrobatmenu#1#2{%
  \@pdfm@mark{%
    bann <<
      /Type /Annot
      /Subtype /Link
      /A <<
        /S /Named
        /N /#1
      >>
      /Border [\@pdfborder]
      /C [\@menubordercolor]
    >>%
   }%
  \Hy@colorlink{\@menucolor}#2\Hy@endcolorlink
  \@pdfm@mark{eann}%
}
\def\@pdfborder{0 0 1}
\def\@menubordercolor{1 0 0}
\def\@menucolor{red}
%
%



\def\marginlabel#1{}

\def\showlabelsabove{
\font\labelfont=cmss10 at 6pt
\def\marginlabel##1{\rlap{\smash{\raise 10pt\hbox{\labelfont##1}}}}
}

\newcount\seccount
\newcount\proccount
\seccount=0
\proccount=0

\def\stdskip{\vskip 9pt plus3pt minus 3pt}
\def\stdbreak{\par\removelastskip\penalty-100\stdskip}

\def\proof{\stdbreak\noindent \colored{\proofcolor}{{\sl Proof. }}}

\def\qed{\vrule height 1.2ex width .9ex depth .1ex}

\def\Box{
  \ifmmode\eqno\qed
  \else\ifvmode\removelastskip\line{\hfil\qed}
  \else\unskip\quad\hskip-\hsize
    \hbox{}\hskip\hsize minus 1em\qed\par
  \fi\stdbreak\fi}

\def\references{
  \removelastskip
  \widow{.05}
  \vskip 24pt plus 6pt minus 6 pt
  \parindent=0pt
  \frenchspacing
  \leftline{\sectionfont \colored{\sectioncolor}{References}}
  \ifhypers\global\advance\seccount by 1\global\proccount=0\relax\edef\numtoks{\number\seccount}%
  \hbox{}\fi
  \nobreak\stdskip\noindent}

\def\ifempty#1#2\endB{\ifx#1\endA}
\def\makeref#1#2#3{\ifempty#1\endA\endB\else\forward{#1}{#2}{#3}\fi}

\outer\def\section#1 #2\par{
  \removelastskip
  \global\advance\seccount by 1
  \global\proccount=0\relax
                \edef\numtoks{\number\seccount}
  \htmlanchor{#1}{\makeref{#1}{Section}{\numtoks}}
  \widow{.05}
  \vskip 24pt plus 6pt minus 6 pt
  \message{#2}
  \ifhypers\hbox{}\fi
  \leftline{\marginlabel{#1}\sectionfont\colored{\sectioncolor}{\numtoks}\quad \colored{\sectioncolor}{#2}}
  \nobreak\stdskip}

\def\proclamation#1#2{
  \outer\expandafter\def\csname#1\endcsname##1 ##2\par{
  \stdbreak
  \global\advance\proccount by 1
  \edef\numtoks{\number\seccount.\number\proccount}
  \htmlanchor{##1}{\makeref{##1}{#2}{\numtoks}}
  \noindent{\marginlabel{##1}\bf \colored{\proclaimcolor}{#2} \colored{\proclaimcolor}{\numtoks}\enspace}
  {\sl##2\par}
  \stdbreak}}

\def\othernumbered#1#2{
  \outer\expandafter\def\csname#1\endcsname##1{
  \stdbreak
  \global\advance\proccount by 1
  \edef\numtoks{\number\seccount.\number\proccount}
  \htmlanchor{##1}{\makeref{##1}{#2}{\numtoks}}
  \noindent{\marginlabel{##1}\bf \colored{\proclaimcolor}{#2} \colored{\proclaimcolor}{\numtoks}\enspace}}}

\proclamation{definition}{Definition}
\proclamation{lemma}{Lemma}
\proclamation{proposition}{Proposition}
\proclamation{theorem}{Theorem}
\proclamation{corollary}{Corollary}
\proclamation{conjecture}{Conjecture}

\othernumbered{example}{Example}
\othernumbered{remark}{Remark}
\othernumbered{construction}{Construction}
\othernumbered{problem}{Problem}

\def\figure#1{
 \global\advance\figcount by 1
 \goodbreak
 \midinsert#1\smallskip
 \centerline{Figure~\number\figcount}
 \endinsert}

\def\capfigure#1#2{
 \global\advance\figcount by 1
 \goodbreak
 \midinsert#1\smallskip
 \vbox{\small\noindent {\bf Figure~\number\figcount:} #2}
 \endinsert}

\def\capfigurepair#1#2#3#4{
 \goodbreak
 \midinsert
 #1\smallskip
 \global\advance\figcount by 1
 \vbox{\small\noindent {\bf Figure~\number\figcount:} #2}
 \vskip 12pt
 #3\smallskip
 \global\advance\figcount by 1
 \vbox{\small\noindent {\bf Figure~\number\figcount:} #4}
 \endinsert}


\def\baretable#1#2{
\vbox{\offinterlineskip\halign{
 \strut\kern #1\hfil##\kern #1
 &&\kern #1\hfil##\kern #1\cr
 #2
}}}

\def\gridtablesub#1#2#3{
\vbox{\offinterlineskip\halign{
 \strut\vrule\kern #1\hfil##\hfil\kern #2\vrule
 &&\kern #1\hfil##\kern #2\vrule\cr
 \noalign{\hrule}
 #3
 \noalign{\hrule}
}}}






\newif\iftextures
\input epsf

\newcount\figcount
\figcount=0
\newcount\figxscale
\newcount\figyscale
\newcount\figxoffset
\newcount\figyoffset
\newbox\drawing
\newcount\drawbp
\newdimen\drawx
\newdimen\drawy
\newdimen\ngap
\newdimen\sgap
\newdimen\wgap
\newdimen\egap

\def\drawbox#1#2#3{\vbox{
  \epsfgetbb{#2.eps} 
  \drawbp=\epsfurx
  \advance\drawbp by-\epsfllx\relax
  \multiply\drawbp by #1
  \divide\drawbp by 100
  \drawx=\drawbp bp
  \drawbp=\epsfury
  \advance\drawbp by-\epsflly\relax
  \multiply\drawbp by #1
  \divide\drawbp by 100
  \drawy=\drawbp bp
  \iftextures
  		\figxscale=#1
    \multiply\figxscale by 10
    \setbox\drawing=\vbox to \drawy{\vfil
      \hbox to \drawx{\special{illustration #2.eps scaled
\number\figxscale}\hfil}}
  \else 
    \figxoffset=-\epsfllx
    \multiply\figxoffset by#1
    \divide\figxoffset by100
    \figyoffset=-\epsflly
    \multiply\figyoffset by#1
    \divide\figyoffset by100
    \setbox\drawing=\vbox to \drawy{\vfil
      \hbox to \drawx{\includegraphics{#2.eps}\hfil}}
  \fi
  \setbox\drawing=\vbox{\offinterlineskip\box\drawing\kern 0pt}
   \ngap=0pt \sgap=0pt \wgap=0pt \egap=0pt
  \setbox0=\vbox{\offinterlineskip
    \box\drawing \ifgridlines\drawgrid\drawx\drawy\fi #3}
  \kern\ngap\hbox{\kern\wgap\box0\kern\egap}\kern\sgap}}

\def\draw#1#2#3{
  \setbox\drawing=\drawbox{#1}{#2}{#3}
  \global\advance\figcount by 1
  \edef\numtoks{\number\figcount}
  \makeref{fig:#2}{Figure}{\numtoks}
  \goodbreak
  \midinsert
  \centerline{\ifgridlines\boxgrid\drawing\fi\box\drawing}
  \smallskip
  \vbox{\offinterlineskip
    \centerline{Figure~\number\figcount}
    \smash{\marginlabel{#2}}}
  \endinsert}

\def\capdraw#1#2#3#4{
  \setbox\drawing=\drawbox{#1}{#2}{#3}
  \global\advance\figcount by 1
  \edef\numtoks{\number\figcount}
  \makeref{fig:#2}{Figure}{\numtoks}
  \goodbreak
  \midinsert
  \centerline{\ifgridlines\boxgrid\drawing\fi\box\drawing}
  \smallskip
  \vbox{\offinterlineskip
    \vskip 4pt
    \vbox{\centerline{\lineskip=3pt\small\noindent
                       {\bf Figure~\number\figcount:} #4}}
    \smash{\marginlabel{fig:#2}}}
  \endinsert}

\def\capdrawpair#1#2#3#4#5#6#7#8{
  \goodbreak
  \midinsert
  \setbox\drawing=\drawbox{#1}{#2}{#3}
  \global\advance\figcount by 1
  \edef\numtoks{\number\figcount}
  \makeref{fig:#2}{Figure}{\numtoks}
  \centerline{\ifgridlines\boxgrid\drawing\fi\box\drawing}
  \smallskip
  \vbox{\offinterlineskip
    \vskip 4pt
    \vbox{\lineskip=3pt\small\noindent {\bf Figure~\number\figcount:} #4}
    \smash{\marginlabel{fig:#2}}}
  \vskip 12pt
  \setbox\drawing=\drawbox{#5}{#6}{#7}
  \global\advance\figcount by 1
  \edef\numtoks{\number\figcount}
  \makeref{fig:#6}{Figure}{\numtoks}
  \centerline{\ifgridlines\boxgrid\drawing\fi\box\drawing}
  \smallskip
  \vbox{\offinterlineskip
    \vskip 4pt
    \vbox{\lineskip=3pt\small\noindent {\bf Figure~\number\figcount:} #8}
    \smash{\marginlabel{fig:#6}}}
  \endinsert}

\def\nextfigtoks{%
  \advance\figcount by 1%
  \edef\numtoks{\number\figcount}%
  \advance\figcount by -1}

\newif\ifgridlines
\newbox\figtbox
\newbox\figgbox
\newdimen\figtx
\newdimen\figty

\newdimen\bwd
\bwd=2sp 

\def\hline#1{\vbox{\smash{\hbox to #1{\leaders\hrule height \bwd\hfil}}}}

\def\vline#1{\hbox to 0pt{%
  \hss\vbox to #1{\leaders\vrule width \bwd\vfil}\hss}}

\def\clap#1{\hbox to 0pt{\hss#1\hss}}
\def\vclap#1{\vbox to 0pt{\offinterlineskip\vss#1\vss}}

\def\hstutter#1#2{\hbox{%
  \setbox0=\hbox{#1}%
  \hbox to #2\wd0{\leaders\box0\hfil}}}

\def\vstutter#1#2{\vbox{
  \setbox0=\vbox{\offinterlineskip #1}
  \dp0=0pt
  \vbox to #2\ht0{\leaders\box0\vfil}}}

\def\crosshairs#1#2{
  \dimen1=.002\drawx
  \dimen2=.002\drawy
  \ifdim\dimen1<\dimen2\dimen3\dimen1\else\dimen3\dimen2\fi
  \setbox1=\vclap{\vline{2\dimen3}}
  \setbox2=\clap{\hline{2\dimen3}}
  \setbox3=\hstutter{\kern\dimen1\box1}{4}
  \setbox4=\vstutter{\kern\dimen2\box2}{4}
  \setbox1=\vclap{\vline{4\dimen3}}
  \setbox2=\clap{\hline{4\dimen3}}
  \setbox5=\clap{\copy1\hstutter{\box3\kern\dimen1\box1}{6}}
  \setbox6=\vclap{\copy2\vstutter{\box4\kern\dimen2\box2}{6}}
  \setbox1=\vbox{\offinterlineskip\box5\box6}
  \smash{\vbox to #2{\hbox to #1{\hss\box1}\vss}}}

\def\boxgrid#1{\rlap{\vbox{\offinterlineskip
  \setbox0=\hline{\wd#1}
  \setbox1=\vline{\ht#1}
  \smash{\vbox to \ht#1{\offinterlineskip\copy0\vfil\box0}}
  \smash{\vbox{\hbox to \wd#1{\copy1\hfil\box1}}}}}}

\def\drawgrid#1#2{\vbox{\offinterlineskip
  \dimen0=\drawx
  \dimen1=\drawy
  \divide\dimen0 by 10
  \divide\dimen1 by 10
  \setbox0=\hline\drawx
  \setbox1=\vline\drawy
  \smash{\vbox{\offinterlineskip
    \copy0\vstutter{\kern\dimen1\box0}{10}}}
  \smash{\hbox{\copy1\hstutter{\kern\dimen0\box1}{10}}}}}

\def\figtext#1#2#3#4#5{
  \setbox\figtbox=\vbox{\hbox{#5}\kern 0pt}
  \figtx=-#3\wd\figtbox \figty=-#4\ht\figtbox
  \advance\figtx by #1\drawx \advance\figty by #2\drawy
  \dimen0=\figtx \advance\dimen0 by\wd\figtbox \advance\dimen0 by-\drawx
  \ifdim\dimen0>\egap\global\egap=\dimen0\fi
  \dimen0=\figty \advance\dimen0 by\ht\figtbox \advance\dimen0 by-\drawy
  \ifdim\dimen0>\ngap\global\ngap=\dimen0\fi
  \dimen0=-\figtx
  \ifdim\dimen0>\wgap\global\wgap=\dimen0\fi
  \dimen0=-\figty
  \ifdim\dimen0>\sgap\global\sgap=\dimen0\fi
  \smash{\rlap{\vbox{\offinterlineskip
    \hbox{\hbox to \figtx{}\ifgridlines\boxgrid\figtbox\fi\box\figtbox}
    \vbox to \figty{}
    \ifgridlines\crosshairs{#1\drawx}{#2\drawy}\fi
    \kern 0pt}}}}


\def\hpad#1#2#3{\hbox{\kern #1\hbox{#3}\kern #2}}
\def\vpad#1#2#3{\setbox0=\hbox{#3}\vbox{\kern #1\box0\kern #2}}




\def\stack#1#2#3{\vbox{\offinterlineskip
  \setbox2=\hbox{#2}
  \setbox3=\hbox{#3}
  \dimen0=\ifdim\wd2>\wd3\wd2\else\wd3\fi
  \hbox to \dimen0{\hss\box2\hss}
  \kern #1
  \hbox to \dimen0{\hss\box3\hss}}}


\def\hexp#1{%
  \setbox0=\hbox{${}^{#1}$}%
  \hbox to .5\wd0{\box0\hss}}

\def\hsub#1{%
  \setbox0=\hbox{${}_{#1}$}%
  \hbox to .5\wd0{\box0\hss}}



\def\bmatrix#1#2{{\left(\vcenter{\halign
  {&\kern#1\hfil$##\mathstrut$\kern#1\cr#2}}\right)}}

\def\rightarrowmat#1#2#3{
  \setbox1=\hbox{\small\kern#2$\bmatrix{#1}{#3}$\kern#2}
  \,\vbox{\offinterlineskip\hbox to\wd1{\hfil\copy1\hfil}
    \kern 3pt\hbox to\wd1{\rightarrowfill}}\,}

\def\leftarrowmat#1#2#3{
  \setbox1=\hbox{\small\kern#2$\bmatrix{#1}{#3}$\kern#2}
  \,\vbox{\offinterlineskip\hbox to\wd1{\hfil\copy1\hfil}
    \kern 3pt\hbox to\wd1{\leftarrowfill}}\,}

\def\rightarrowbox#1#2{
  \setbox1=\hbox{\kern#1\hbox{\small #2}\kern#1}
  \,\vbox{\offinterlineskip\hbox to\wd1{\hfil\copy1\hfil}
    \kern 3pt\hbox to\wd1{\rightarrowfill}}\,}

\def\leftarrowbox#1#2{
  \setbox1=\hbox{\kern#1\hbox{\small #2}\kern#1}
  \,\vbox{\offinterlineskip\hbox to\wd1{\hfil\copy1\hfil}
    \kern 3pt\hbox to\wd1{\leftarrowfill}}\,}








\def\quiremacro#1#2#3#4#5#6#7#8#9{
  \expandafter\def\csname#1\endcsname##1{
  \ifnum\srcdepth=1
  \magnification=#2
  \input quire
  \hsize=#3
  \vsize=#4
  \htotal=#5
  \vtotal=#6
  \shstaplewidth=#7
  \shstaplelength=#8
  \hoffset=\htotal
  \advance\hoffset by -\hsize
  \divide\hoffset by 2
  \ifnum\vsize<\vtotal
    \voffset=\vtotal
    \advance\voffset by -\vsize
    \divide\voffset by 2
  \fi
  \advance\voffset by #9
  \shhtotal=2\htotal
  \baselineskip=13pt
  \emergencystretch = 0.05\hsize
  \horigin=0.0truein
  \vorigin=0.0truein
  \shthickness=0pt
  \shoutline=0pt
  \shcrop=0pt
  \shvoffset=-1.0truein
  \ifnum##1>0\quire{#1}\else\qtwopages\fi
  \fi
}}



\quiremacro{letterbooklet} 
{1000}{4.79452truein}{7truein}{5.5truein}{8.5truein}{0.01pt}{0.66truein}
{-.0625truein}

\quiremacro{Afourbooklet}
{1095}{5.25truein}{7truein}{421truept}{595truept}{0.01pt}{0.66truein}
{-.0625truein}

\quiremacro{legalbooklet}
{1095}{5.25truein}{7truein}{7.0truein}{8.5truein}{0.01pt}{0.66truein}
{-.0625truein}

\quiremacro{twoupsub} 
{895}{4.5truein}{7truein}{5.5truein}{8.5truein}{0pt}{0pt}{.0625truein}


\quiremacro{Afourviewsub} 
{1000}{5.0228311in}{7.7625571in}{421truept}{595truept}{0.1pt}{0.5\vtotal}
{-.0625truein}


\quiremacro{viewsub}
{1095}{5.5truein}{8.5truein}{461truept}{666truept}{0.1pt}{0.5\vtotal}
{-.125truein}


\newcount\countA
\newcount\countB
\newcount\countC

\def\monthname{\begingroup
  \ifcase\number\month
    \or January\or February\or March\or April\or May\or June\or
    July\or August\or September\or October\or November\or December\fi
\endgroup}

\def\dayname{\begingroup
  \countA=\number\day
  \countB=\number\year
  \advance\countA by 0 
  \advance\countA by \ifcase\month\or
    0\or 31\or 59\or 90\or 120\or 151\or
    181\or 212\or 243\or 273\or 304\or 334\fi
  \advance\countB by -1995
  \multiply\countB by 365
  \advance\countA by \countB
  \countB=\countA
  \divide\countB by 7
  \multiply\countB by 7
  \advance\countA by -\countB
  \advance\countA by 1
  \ifcase\countA\or Sunday\or Monday\or Tuesday\or Wednesday\or
    Thursday\or Friday\or Saturday\fi
\endgroup}

\def\timename{\begingroup
   \countA = \time
   \divide\countA by 60
   \countB = \countA
   \countC = \time
   \multiply\countA by 60
   \advance\countC by -\countA
   \ifnum\countC<10\toks1={0}\else\toks1={}\fi
   \ifnum\countB<12 \toks0={\sevenrm AM}
     \else\toks0={\sevenrm PM}\advance\countB by -12\fi
   \relax\ifnum\countB=0\countB=12\fi
   \hbox{\the\countB:\the\toks1 \the\countC \thinspace \the\toks0}
\endgroup}

\def\timestamp{\dayname, \the\day\ \monthname\ \the\year, \timename}


\print



\def\COMMENT#1\par{\bigskip\hrule\smallskip#1\smallskip\hrule\bigskip}

\def\enma#1{{\ifmmode#1\else$#1$\fi}}

\def\mathbb#1{{\bbold #1}}
\def\mathbf#1{{\bf #1}}


\def\PP{\enma{\mathbb{P}}}
\def\HH{\enma{\mathbb{H}}}



%
\def\boldN{\enma{\mathbf{N}}}
\font\tengoth=eufm10  \font\fivegoth=eufm5
\font\sevengoth=eufm7
\newfam\gothfam  \scriptscriptfont\gothfam=\fivegoth 
\textfont\gothfam=\tengoth \scriptfont\gothfam=\sevengoth

\def\set#1{\enma{\{#1\}}}


\def\codim{\mathop{\rm cod}\nolimits}

\def\lspan{{\mathop{\rm span}\nolimits}}
\def\reg{\mathop{\rm reg}\nolimits}
\def\Im{\mathop{\rm Im}\nolimits}

\def\ker{\mathop{\rm ker}\nolimits}
\def\deg{\mathop{\rm deg}\nolimits}

\def\dim{\mathop{\rm dim}\nolimits}

%

%

%
%

\newsymbol\boxtimes1202


\input diagrams.tex
\proclamation{question}{Question}

\def\codim{\mathop{\rm codim}\nolimits}

\def\lspan{{\mathop{\rm span}\nolimits}}
\def\reg{\mathop{\rm reg}\nolimits}
\def\Im{\mathop{\rm Im}\nolimits}

\def\ker{\mathop{\rm ker}\nolimits}
\def\deg{\mathop{\rm deg}\nolimits}

\def\dim{\mathop{\rm dim}\nolimits}

\def\red{\enma{\rm red}}


%



\def\enma#1{{\ifmmode#1\else$#1$\fi}}
\def\mathbb#1{{\bbold #1}}
\def\mathbf#1{{\bf #1}}


\def\PP{\enma{\mathbb{P}}}
\def\HH{\enma{\mathbf{H}}}


\let\I\cII

\let\O\cOO

%
\def\boldN{\enma{\mathbf{N}}}
\font\tengoth=eufm10  \font\fivegoth=eufm5
\font\sevengoth=eufm7
\newfam\gothfam  \scriptscriptfont\gothfam=\fivegoth 
\textfont\gothfam=\tengoth \scriptfont\gothfam=\sevengoth

\def\set#1{\enma{\{#1\}}}

\forward{2-regular schemes}{Section}{1}
\forward{small and geomsmall}{Section}{2}
\forward{LJ section}{Section}{3}
\forward{basic properties}{Section}{4}
\forward{chordal}{Section}{5}
\forward{equations}{Section}{6}

\forward {quartic in P4}{Example}{0.5}
\forward{project 2-regular}{Corollary}{0.7}
\forward{2-reg reduced 2-reg}{Corollary}{0.8}

\forward{2-regular gives small}{Theorem}{1.2}

\forward{small basics}{Proposition}{2.1}
\forward{small char}{Theorem}{2.2}
\forward{subtracting disjoint set}{Corollary}{2.3}

\forward{linearly joined equivalences}{Proposition}{3.1}
\forward{LJ for schemes and spans}{Proposition}{3.4}

\forward{oberwolfach 2-regular}{Theorem}{4.7}
\forward{small projects small}{Proposition}{4.4}
\forward{iso projection and plane section}{Proposition}{4.8}
\forward{2-regular secants}{Proposition}{4.10}
\forward{not small}{Corollary}{4.11}

\forward{weights of forests}{Theorem}{5.1}

\forward{lin join equations}{Theorem}{6.1}
\bigskip
\centerline{\titlefont Small Schemes and Varieties of Minimal Degree}
\medskip
\centerline {by}
\smallskip
\centerline {\bf D. Eisenbud, M. Green, K. Hulek, and S. Popescu}
\footnote{}{\rm The authors
are grateful to BIRS, IPAM, MSRI, the DFG and the NSF for
hospitality and support during the preparation of this work.}

\bigskip\bigskip

\noindent {\bf Abstract:} 
{\narrower
\small 
We prove that if $X\subset \PP^r$ is any 2-regular scheme (in the
sense of Castelnuovo-Mumford) then $X$ is {\it small\/}. This means
that if $L$ is a linear space and $Y:= L\cap X$ is finite, then $Y$ is
{\it linearly independent\/} in the sense that the dimension of the
linear span of $Y$ is $1+\deg Y$. The converse is true and well-known
for finite schemes, but false in general.  The main result of this
paper is that the converse, ``small implies 2-regular'', is also true
for reduced schemes (algebraic sets). This is proven by means of a
delicate geometric analysis, leading to a complete classification: 
we show that the components of a small algebraic set
are varieties of minimal degree, meeting in a particularly
simple way.  From the classification one can show 
that if $X\subset \PP^r$ is 2-regular,
then so is $X_\red$, and so also is the projection of
$X$ from any point of $X$.

Our results extend the Del Pezzo-Bertini classification of varieties
of minimal degree, the characterization of these as the varieties of
regularity 2 by Eisenbud-Goto, and the construction of 2-regular
square-free monomial ideals by Fr\"oberg.

}

\bigskip\bigskip\noindent

Throughout this paper we will work with projective schemes $X\subset \PP^r$
over an algebraically closed field $k$. The (Castelnuovo-Mumford) regularity
of $X\subset \PP^r$ is a basic homological measure of the complexity of $X$
and its embedding in $\PP^r$ that gives a bound for the degrees of the generators
of the defining ideal $I_X$ of $X$ and for many other invariants.
The only schemes of regularity 1 are the linear spaces; but no
classification is known for projective schemes of regularity 2. 

In this paper we prove a structure 
theorem for reduced 2-regular
schemes, showing that their irreducible components are
varieties of minimal degree and characterizing how
these components can meet. We also show
that the reduced structure on any 2-regular scheme is
2-regular, and thus we obtain a complete  description of the
reduced structures on 2-regular schemes. (Since a
high Veronese re-embedding of any zero-dimensional
scheme is 2-regular, one cannot
hope to characterize the isomorphism types of all 2-regular
non-reduced schemes.)

Before stating our results we review some basic notions.
For any subscheme $X\subset \PP^r$ we write $\lspan(X)$ for the smallest
linear subspace of $\PP^r$ containing $X$.  Recall that every variety
($\equiv$ reduced irreducible scheme)  $X\subset
\PP^r$ satisfies the condition 
$$
\deg(X)\ge 1+\codim(X, \lspan(X))\leqno{(*)}
$$ 
(see for instance Mumford [1976, Corollary 5.13]). We say that the variety
$X\subset\PP^r$ has {\it minimal degree\/} (more precisely,
minimal degree in its span) if equality holds. 
Surfaces of
minimal degree were classified by Del Pezzo [1886], and the
classification was extended to all dimensions by Bertini
[1907] (see Eisenbud-Harris [1987] for a modern account):  

\theorem{del Pezzo-Bertini} A projective variety of minimal degree 
in its span is either a linear space, a quadric hypersurface
in a linear space, a rational
normal scroll, or a cone over the Veronese surface in $\PP^5$.

This classification was extended to equidimensional algebraic sets
that are connected in codimension 1 -- the ones for which ``minimal
degree'' is a good notion -- by Xamb\'o [1981]. For more general
algebraic sets  it is not clear that there should exist any interesting
generalization of the equality in $(*)$ above. Nevertheless,
the notion of smallness
is just such a generalization.

Varieties of minimal degree were characterized cohomologically by
Eisenbud-Goto [1984], and their result offers a
way to generalize the hypothesis of the Del Pezzo-Bertini
Theorem to all projective schemes. To state their result recall that 
$X\subset \PP^r$ is said to have {\it regularity $d$}, or to 
be {\it $d$-regular,\/} in the sense of
Castelnuovo-Mumford, if the ideal sheaf $\I_X$ satisfies
$H^i(\I_X(d-i))=0$ for all $i>0$, or equivalently, if the $j$-th
syzygies of the homogeneous ideal $I_X$ are generated in degrees 
$\leq d+j$ for all $j\geq 0$ (see Eisenbud-Goto [1984], or Eisenbud [2004]
for a proof of the equivalence).

\theorem{homol char of min deg}
A variety $X\subset \PP^r$ has minimal degree in its linear span if
and only if $X$ is $2$-regular.

An old argument of Lazarsfeld (see for instance [2004]), recently
refined by Sidman [2002], Caviglia [2003], Eisenbud-Green-Hulek-Popescu [2004] 
and others, shows that if $X\subset\PP^r$ is $d$-regular and $\Lambda$ is a linear 
subspace such that $X\cap \Lambda$ has dimension 0, then $X\cap \Lambda$ is 
also $d$-regular.  When $d=2$ we can rephrase this geometrically:

We say that a finite scheme $Y\subset \PP^r$ is {\it linearly independent\/} if the dimension of
the linear span of $Y$ is $1+\deg Y$.  We say that a scheme $X\subset
\PP^r$ is {\it small\/} if, for every linear subspace $\Lambda\subset \PP^r$
such that $Y= \Lambda\cap X$ is finite, the scheme $Y$ is linearly
independent. (An alternative definition of smallness by a more general
property of intersections is given in \ref{small char}.) 
Lazarsfeld's argument gives:

\proposition{smallness of 2-regular} Any 2-regular scheme
$X\subset \PP^r$ is small.

Our main results are that the converse holds in the reduced case,
and that small reduced schemes
have a simple inductive classification.
 To state the classification,
we say that  a sequence of closed subschemes  $X_1,\dots, X_n\subset \PP^r$ is 
{\it linearly joined\/} if, for all $i=1,\dots,n-1$, we have
$$
(X_1\cup\ldots\cup X_i)\cap X_{i+1}
=
\lspan(X_1\cup\ldots\cup X_i)\cap \lspan(X_{i+1}).
$$

\theorem{main} Let $X\subset \PP^r$ be an algebraic
set. The following conditions are equivalent:
\item{$(a)$} $X$ is small.
\item{$(b)$} $X$ is $2$-regular.
\item{$(c)$} $X=X_1\cup\ldots\cup X_n$, where 
$X_1,\dots,X_n$ is a linearly joined sequence of 
varieties of minimal degree.

The implication $(c)\Rightarrow(b)$ is easy (see \ref{linearly joined
equivalences}) while $(b)\Rightarrow (a)$ is a special case of
\ref{smallness of 2-regular} (see also \ref{2-regular schemes} for details).  
Most of this paper is occupied with the proof of the implication
$(a)\Rightarrow(c)$, which requires a delicate geometric analysis of
the notion of smallness in the style of classical projective 
geometry. One of the things that makes the argument subtle
is the fact that the linearly joined property of a sequence
of varieties is strongly
dependent on the ordering, as the following example shows.

\example{quartic in P4} Let $L_0$ be a line in $\PP^4$, and let
$L_1,L_2,L_3$ be 3 general lines that meet $L_0$. The union 
$X=\bigcup_{i=0}^3L_i$ is $2$-regular; in fact
it is connected in codimension 1, has minimal degree, and is
a degeneration of a rational normal quartic curve. 
As required by \ref{main}, $X$
can be written as the union of a linearly joined sequence of varieties
$L_0, L_1, L_2, L_3$, which are trivially
of minimal degree in their spans. 

On the other hand, the reverse sequence $L_3, L_2, L_1, L_0$
is not linearly joined. Indeed,
the subset $Y=L_3\cup L_2\cup L_1$
is not $2$-regular: since
$Y$ meets the line $L_0$ in three points,
the ideal of $Y$ requires a cubic generator. It
is easy to check that there is
no enumeration of the components of $X$ as a linearly joined sequence
for which the reverse sequence is linearly joined.
\medskip

In the special case where $X$ is a union of
 {\it coordinate\/} subspaces,
the equivalence of parts $(b)$ and $(c)$ of \ref{main}
had been proved by Fr\"oberg [1985, 1988] as an
application of Stanley-Reisner
theory. Unions of coordinate spaces correspond to simplicial complexes.
Using earlier results of Dirac [1961] and Fulkerson-Gross [1965],
Fr\"oberg showed that a simplicial complex corresponds to
a $2$-regular set if and only if it is the clique complex of
a chordal graph. (We reprove this and give a generalization
in Eisenbud-Green-Hulek-Popescu [2004]. See also Herzog, Hibi, and 
Zheng [2003] for a related path to Dirac's theorem.) 
The orderings described in part $(c)$ of \ref{main}
are called {\it perfect elimination orderings\/}
in this context. See Blair-Peyton [1993] for a survey.

Properties that are easy to check for algebraic sets satisfying one of
the conditions of \ref{main} may be quite obscure for sets satisfying
another. \ref{main} has a number of surprising algebraic and geometric
consequences based on this observation:

\corollary{2 components} If $X\subset \PP^r$ is a $2$-regular algebraic
set, then the union of any two irreducible components of $X$ is again $2$-regular.

By \ref{quartic in P4} the same cannot be said of a union of three components.

\proof By \ref{main}, the irreducible components
of $X$ are of minimal degree in their spans. 
and thus also $2$-regular by \ref{homol char of min deg}.
From 
\ref{linearly joined equivalences}
we see that any 2 components are again linearly joined.
The result follows by applying \ref{main} once again.\Box

\corollary{project 2-regular} Let $X\subset\PP^r$
be a $2$-regular algebraic set. If $p\in X$ is a point and
$\pi_p$ denotes the linear projection from $p$, then
$\pi_p(X)\subset\PP^{r-1}$ is $2$-regular.

By \ref{main} a similar statement holds with
``small'' in place of ``$2$-regular''.

\proof If 
$X_1,\dots,X_n$ is a linearly joined sequence of varieties of minimal degree,
then $\pi_p(X_1),\dots,\pi_p(X_n)$ is a linearly joined sequence 
of varieties of minimal degree. \ref{main} completes the proof.\Box

\corollary{2-reg reduced 2-reg} If $X\subset\PP^r$ is $2$-regular, 
then $X_{\rm{red}}\subset\PP^r$ is also
$2$-regular.

\proof By Eisenbud-Green-Hulek-Popescu [2004, Theorem 1.6] (see
also \ref{2-regular gives small}) the scheme $X\subset\PP^r$  is
small. It follows that $X_{\rm{red}}$ is also small 
(\ref{small basics}). By \ref{main}, $X_\red$ is $2$-regular.\Box

\corollary{oberwolfach-small} 
Let $X_1,\dots, X_n\subset \PP^r$ be a collection
of varieties of minimal degree in their spans.
 The
union $\bigcup_iX_i$ is small if and only if each
pair $X_i, X_j$ is linearly joined and the 
union of the linear spans $\bigcup_i \lspan(X_i)$ is small.

Of course a similar statement will hold for $2$-regularity in place
of smallness. That version is actually one of the key ingredients in the
proof of \ref{main}.
 
\proof Use \ref{main} and \ref{LJ for schemes and spans}\Box

The plan of the paper is as follows. In Sections 
\refn{2-regular schemes}, \refn{small and geomsmall} and
\refn{LJ section} we establish basic properties of $2$-regular sets,
small projective schemes, and linearly joined sequences of projective
schemes. Of particular interest is the B\'ezout type
result, \ref{small char}: 
If $X\subset \PP^r$ is a small subscheme and $\Lambda\subset \PP^r$ is any
linear space, then the sum of the degrees of the irreducible
components (reduced or not, but not embedded) of $X\cap \Lambda$ is
bounded by $\codim(X\cap\Lambda,\Lambda)+1$. The 
results in these sections are necessary for the proof of \ref{main},
which is carried out in \ref{basic properties}.

The argument of Lazarsfeld showing that plane sections
of $X\cap L$ of a $2$-regular scheme $X$ are small works even if $X$ is
not $2$-regular, but only has an ideal generated by
quadrics having only linear syzygies for at least $\dim L$
steps. With slightly stronger hypotheses one can prove a little
more; for example that that the syzygies of $X\cap L$ come from 
syzygies of $X$  by restriction. See \ref{2-regular gives small}, and
Eisenbud-Green-Hulek-Popescu [2004]. 

It follows from \ref{oberwolfach-small} 
that the condition that
an algebraic set $X\subset \PP^r$ be small (or $2$-regular) 
has a ``local'' part, that the $X_i$ of $X$
should be of minimal degree in their spans, and pairwise linearly joined;
and a ``global'', or combinatorial part, that
the subspace arrangement $\bigcup_i \lspan(X_i)$ be small.
In \ref{chordal} we will study 
small subspace arrangements. In particular, we describe
the orderings of subspaces that make them into a 
linearly joined sequence in terms of certain spanning forests
of the intersection graph of $Y$.

Finally, in \ref{equations} we discuss how to find
generators for the ideal of a reduced
$2$-regular projective scheme. In particular, 
we prove that any $2$-regular union of linear spaces
has ideal generated by products of linear forms. The proof of
\ref{lin join equations} also provides a free
resolution for the ideal of a $2$-regular
algebraic set, and one can ensure that the first
two terms of the resolution are minimal, obtaining 
a formula for the number of generators of the ideal.
Some results along this line have also been obtained
by Barile and Morales [2000, 2003].

\smallskip
We thank Aldo Conca, Harm Derksen, Mark Haiman, J\"urgen Herzog, 
Bill Oxbury, Kristian Ranestad, Jerzy Weyman, and Sergey Yuzvinsky 
for useful discussions. Behind the scenes, the program 
{\it Macaulay2\/} of Dan Grayson and Mike Stillman
has been extremely useful to us in trying to understand the 
relations between the linear syzygies and geometry.

\section{2-regular schemes} 2-regular schemes

The regularity of the ideal sheaf of a closed subscheme
$X\subset\PP^r$ is $\leq 1$ if and only if $X$ is defined by linear
forms, so $X$ is a linear subspace in this case.

By contrast, {\it any\/} finite scheme can be embedded as a
$2$-regular scheme. In fact, we see from the definition of regularity
that a high Veronese re-embedding of any given embedding of a
zero-dimensional scheme is $2$-regular.  
In fact a complete characterization of $2$-regular embeddings
of zero-dimensional schemes is well-known:

\proposition{2-reg 0-dim}
A zero-dimensional nondegenerate scheme $X\subset \PP^r$ is
$2$-regular if and only if $\deg X = 1+\dim \lspan(X)$.

\proof The cohomological condition for regularity is equivalent to saying
that $X$ imposes $\deg X$ independent conditions on linear forms,
that is, the span of $X$ has dimension $\deg X-1$.
\Box

We will often use the fact that a zero-dimensional
plane section of a $2$-regular scheme is again $2$-regular.
As described in the introduction, this
follows from an argument of Lazarsfeld, 
and many others have studied it recently.
Here is a version sufficient for our purposes in this
paper.
We say that an ideal $I\subset S$ has
{\it 2-linear resolution for at least $p$ steps\/}
if the $i$-th syzygies of $I$ are
generated in degrees $\leq i+1$ for $i=0,\dots,p-1$.
For example, if $I$ contains no linear forms, 
this means that the  minimal free resolution of I has the form
$$
\cdots \oplus S(-i)^{\beta_{p, i}}\rTo S(-p-1)^{\beta_{p-1}}\rTo \cdots
\rTo \oplus S(-2)^{\beta_0}\rTo I \rTo 0.
$$
We also say in this case that $I$ satisfies property $\boldN_{2,p}$;
this is the terminology used in Eisenbud-Green-Hulek-Popescu [2004].

\theorem{2-regular gives small} Let $S$ be the polynomial ring
on $r+1$ variables, and suppose that $X\subset\PP^r$ has homogeneous
ideal $I_X\subset S$ such that $I_X$ is generated by quadrics and has
linear resolution for at least $p$ steps.  If $\Lambda\subset\PP^r$ is
a linear subspace with $\codim(\Lambda\cap X, \lspan(\Lambda\cap
X))\le p-1$, then $\Lambda\cap X$ is $2$-regular. In particular, any
zero-dimensional plane section of a $2$-regular scheme is $2$-regular.\Box

For an explicit proof see for example
Eisenbud-Green-Hulek-Popescu [2004, Th 1.1]. 
Combining \ref{2-reg 0-dim} and \ref{2-regular gives small} we get as
a corollary the \ref{smallness of 2-regular} in the introduction:

\corollary{2-regular gives small cor}
Any $2$-regular closed scheme $X\subset \PP^r$ is small.\Box

There is a geometric characterization of reduced $2$-regular schemes of
higher dimension in the Cohen-Macaulay case. By a result of Hartshorne [1962] 
(see Eisenbud [1995]), connectedness in codimension 1 is a
necessary condition for Cohen-Macaulayness. It turns out that for
reduced $2$-regular algebraic sets they are equivalent.

\theorem{2-reg of min deg} Let $X\subset \PP^r$ be an equidimensional
projective scheme, and let $L$ be the linear span of $X$. Suppose that 
$X$ is reduced,and connected in codimension $1$. 
\item{$(a)$} $\deg X\geq \codim(X, L)+1$.
\item{$(b)$} $X$ is $2$-regular if and only if
$\deg X= \codim(X, L)+1$.
\item{$(c)$} If the equivalent conditions in part $(b)$ hold,
then the homogeneous coordinate ring of $X$ is  Cohen-Macaulay. \Box

\proof Part $(a)$ of \ref{2-reg of min deg}
is elementary: the proof works as in the irreducible case
(see Hartshorne [1962]) using the connectedness hypothesis
to guarantee that the plane section is nondegenerate. The rest is
proven in Eisenbud-Goto [1984].\Box

The following Corollary is \ref{main} in the special case of irreducible or
connected-in-codimension $1$ algebraic sets:

\corollary{2-reg of geom small irred} 
Let $X\subset \PP^r$ be a projective scheme, and let $L$ be the linear
span of $X$. Suppose that $X$ is reduced, and connected in
codimension $1$. If $X$ is small then $X$ is $2$-regular. \Box

A well-known regularity criterion that
is essentially due to Mumford [1966, Lecture 14]
will play a central role. (See also the last section of Eisenbud 
[2004] for details.)

\theorem{2-reg criterion} A closed subscheme $X\subset \PP^r$ is $2$-regular
if and only if
\item {$(a)$} For some (respectively, any) hyperplane $H\subset \PP^r$,
defined by a linear form that is locally a nonzerodivisor on $\O_X$,
the scheme $X\cap H$ is $2$-regular, and
\item {$(b)$} $X$ is linearly normal; equivalently, the
restriction map \hfill\break $H^0(\O_{\PP^r}(1))\rTo H^0(\O_X(1))$
is surjective. \Box

\medskip
We will say that a scheme $X\subset \PP^r$ is the direct sum of 
schemes $Z_i\subset \PP^r$ if, for each $i$,
$$
\lspan(Z_i)\cap \lspan (\bigcup_{j\neq i} Z_j)=\emptyset.
$$
In other words the underlying vector space of the span of $X$ is
the direct sum of the underlying vector spaces of the spans
of the $Z_i$. If this holds, then
most questions about $X\subset \PP^r$ can be reduced 
to questions about the $Z_i\subset \lspan(Z_i)$; for instance 
the next result will allow us to 
reduce questions about $2$-regularity to the connected case.

\proposition{direct sums 2-reg} $X\subset \PP^r$ is $2$-regular
if and only if the connected components of $X$ are $2$-regular and
$X$ is the direct sum of its connected components.

\proof The vanishing $H^i(\I_X(2-i))\cong H^{i-1}(\O_X(2-i))=0$,
for all $i\ge 2$, is equivalent to $H^{i-1}(\O_Z(2-i))=0$, for all $i\ge 2$, and
all connected  components $Z$ of $X$. On the other hand we may assume that $X$ 
is nondegenerate. Then if $X=Z_1\cup Z_2$  is a disjoint union, 
the cohomology of the short exact sequence
$$
0\rTo\I_X(1)\rTo\I_{Z_1}(1)\oplus\I_{Z_2}(1)\rTo\O_{\PP^r}(1)\rTo 0
$$
yields that $H^1(\I_X(1))=0$ if and only if both $H^1(\I_{Z_i}(1))=0$, for $i=1,2$,
and $H^0(\I_{Z_1}(1))\oplus H^0(\I_{Z_2}(1))\rTo H^0(\O_{\PP^r}(1))$ is
an isomorphism.\Box

\section{small and geomsmall} Basic properties of small algebraic sets

The following remarks will be used frequently in the sequel. They follow
at once from the definitions.

\proposition{small basics} Let $X\subset \PP^r$ be a small
subscheme.
\item{$(a)$} Any plane section of $X$ is small.
\item{$(b)$} $X_{\red}$ is small.

We could have defined smallness by a more general property of
intersections.  Recall that the {\it geometric degree} of a closed
subscheme $Y\subset\PP^r$ is defined as the sum of the degrees of the
isolated irreducible (not necessarily reduced) components of $Y$.

\theorem{small char} If $X\subset \PP^r$  is a small scheme, and
$L\subset \PP^r$ is a linear space, then the geometric degree of
$X\cap L$ is at most $1+\codim(X\cap L, L)$.

Does \ref{small char} hold with the {\it arithmetic
degree} (sum of all the degrees of isolated {\it and\/} embedded
components) in place of the geometric degree?  See Bayer-Mumford [1993], 
Sturmfels-Trung-Vogel [1995], or Miyazaki-Vogel-Yanagawa [1997] for
the definition and basic properties of such degrees.

\proof By \ref{small basics} we can assume $L=\PP^r$.  
We do induction on $r$, the case $r=0$ being obvious.  Let $Y$ be the
union of the zero-dimensional components of $X$, so that $X=Y\cup Z$
is the disjoint union of $Y$ and a scheme $Z$ whose components have
positive dimension.  If $Y$ spans $\PP^r$, then successively factoring
out elements of the socle of $\O_Y$, we see that $Y$ contains a
subscheme $Y'$ of length $r$ spanning an $(r-1)$-plane $\Lambda\subset
\PP^r$. If $Z$ is non-empty then this plane meets $Z$ nontrivially,
and thus the sum of the degrees of the components of $X\cap \Lambda$
is $> r=\dim \Lambda+1\geq\codim(\Lambda\cap X, \Lambda)$.  
Since $X\cap \Lambda$ is small, this contradicts the inductive hypothesis 
and shows that $Z=\emptyset$. Therefore $X=Y$ is finite, and the desired 
conclusion follows from the definition of smallness.

On the other hand, suppose that $Y$ does not span $\PP^r$.
If $\Lambda$ is a general hyperplane containing $Y$, then
$\Lambda$ meets all the irreducible components $Z_i$ 
of $Z$ in  schemes $Z_i\cap X$ of the same degree as
$Z_i$. Further, the intersection $Z_i\cap Z_j$
of distinct components has dimension strictly less than
that of $Z_i$ or $Z_j$, so the schemes $Z_i\cap Z_j\cap \Lambda$
do not contain any component of $Z\cap \Lambda$. It follows
that the sum of the degrees of the components of $Z\cap \Lambda$
is equal to the corresponding sum for $Z$. By induction on $r$,
the degree of $Y\subset \Lambda$ plus the sum of the degrees
of the components of $Z\cap \Lambda$ is bounded by
$1+\codim(X\cap \Lambda, \Lambda)=1+\codim(X,\PP^r)$.
Thus the same bound holds for $X$. \Box

We isolate a consequence of \ref{small char} for use in the 
proof of the \ref{main}:

\corollary {subtracting disjoint set} Suppose that
 $Y, Z\subset \PP^s$ are disjoint algebraic sets.  
If $Y\cup Z$ is small then both $Y$ and $Z$ are small.

\proof By symmetry it suffices to show that $Y$ is small.
Suppose that $L$ is a linear space that meets $Y$ in a finite
scheme. Since $L\cap(Y\cup Z)$ is the disjoint union of $L\cap Y$
and $L\cap Z$, the geometric degree of $L\cap(Y\cup Z)$ is
at least as great as the length of $L\cap Y$, so we are done
by \ref{small char}.\Box

The conclusion of \ref{small char} may be interpreted as a B\'ezout
type theorem in the case of small varieties.  A special case of a
result of Lazarsfeld says that if $X\subset \PP^r$ is a nondegenerate
subvariety and $\Lambda\subset \PP^r$ is a linear subspace, then the
geometric degree of $X\cap \Lambda$ is bounded above by
$\deg(X)-\codim(X,\PP^r)+\codim(X\cap \Lambda, \Lambda)$ (see Fulton
[1984], Example 12.3.5, which also states Lazarsfeld's more general
result, or Fulton-Lazarsfeld [1982]).  In the case when $X$ is a 
variety of minimal degree this yields \ref{small char}.

\smallskip

Next we analyze the irreducible components, and the relative
positions of pairs of irreducible components, of small
sets. By \ref{small basics} the same results would apply to 
the reduced irreducible components of any small projective scheme.

\proposition{small basic props} 
Let $X\subset\PP^r$ be a small algebraic set.
\item{$(a)$} Any irreducible component $X_i$ of $X$ 
is a variety of minimal degree in its linear span; that is 
$\deg X_i = \dim\lspan(X_i)-\dim X_i+1$.
\item{$(b)$} Any two irreducible components $X_i$ and $X_j$ of $X$ 
are linearly joined; that is,
$
X_i\cap X_j=\lspan(X_i)\cap\lspan(X_j).
$

\proof  Write $X=\bigcup_i X_i$, where $X_i$ are the irreducible components of
$X$, and let $L_i=\lspan(X_i)$ denote the linear span of $X_i$. Any
variety $X_i$ in projective space satisfies $\deg X_i\geq
\dim\lspan(X_i)-\dim X_i+1$, but in our case, $X_i$ is a component of
$L_i\cap X$, so the opposite inequality and thus part $(a)$ follows 
by applying \ref{small char} to the linear space $L_i$.

For part $(b)$, we first apply \ref{small char} to the linear
space $L=\lspan(L_i\cup L_j)$. By construction $X_i$ and $X_j$
are components of $L\cap X$, so the geometric degree
of $L\cap X$ is at least $\deg X_i+\deg X_j$, and so \ref{small char} yields
$$
\deg X_i+\deg X_j \leq 
\dim L_i+\dim L_j-\dim(L_i\cap L_j)-\max(\dim X_i,\dim X_j)+1.
$$
On the other hand, by part $(a)$, $\deg X_i=\dim L_i-\dim X_i+1$
and we deduce
$\dim(L_i\cap L_j) \leq \min (\dim X_i, \dim X_j)-1$.

Now suppose that $X_i\cap X_j\neq L_i\cap L_j$.  If $\dim L_i\cap
L_j>0$, then to simplify we may cut all the schemes concerned by a
general hyperplane $H$.  Because $X_i$ is reduced, irreducible and of
dimension $>1$, the hyperplane section $H\cap X_i$ is again reduced
and irreducible (by Bertini's Theorem) and spans $L_i\cap H$. The same
holds for $X_j$.  Continuing to take hyperplane sections, we may
reduce to the case where the linear space $L_i\cap L_j$ is just 
a point $p$.

If now both $X_i$ and $X_j$ contain $p$, then 
$X_i\cap X_j = \{p\} = L_i\cap L_j$ as desired.
If however $X_i$ does not contain $p$ but
$X_j$ does, then $L_i\cap X$ contains a component
through $p$ and thus 
$$
\deg (X\cap L_i) > \deg X_i
\geq  \dim L_i - \dim X_i +1 \geq \dim L_i - \dim (X\cap L_i)+1,
$$
contradicting \ref{small char}. By symmetry, we may therefore 
assume that neither $X_i$ nor $X_j$ meets $L_i\cap L_j=\{p\}$, 
and we must derive a contradiction.

Let $\Lambda_k\subset L_k$, for $k=i,j$, be general planes containing
$p$ and having $\dim \Lambda_k =\codim(X_k, L_k)$. The scheme
$\Lambda_k\cap X_k$ is zero-dimensional.
Set $\Lambda=\lspan(\Lambda_i\cup\Lambda_j)$. We have
$\Lambda\cap L_i = \Lambda_i+(\Lambda_j\cap L_i) = \Lambda_i$,
and similarly for $\Lambda\cap L_j$, so the components
of $\Lambda_k\cap X_k$, $k=i,j$, are also components of $\Lambda\cap X$.
Thus the geometric degree of $\Lambda\cap X$ is at least
$\deg X_i+\deg X_j$. 

On the other hand, since $\Lambda_i$ and $\Lambda_j$ meet in a point,
the dimension of $\Lambda$ is $\dim \Lambda_i+\dim
\Lambda_j=\codim(X_i, L_i)+\codim(X_j,L_j)=\deg X_i +\deg X_j -2$, 
by part $(a)$. But then \ref{small char} yields that the geometric
degree of $\Lambda\cap X$ is at most $\dim \Lambda+1 = \deg X_i+\deg X_j-1$,
the desired contradiction. This concludes the proof of $(b)$.\Box

\section{LJ section} Linearly joined sequences of schemes

A first connection of the notions of smallness and $2$-regularity with
the notion of linearly joined sequence is provided by the following
observations. A linearly joined sequence
of schemes $X,Y\subset \PP^r$ is linearly joined in either order,
so we simply say that the pair is linearly joined.

\proposition{linearly joined equivalences}
Suppose $X_1,X_2\subset \PP^r$ are linearly joined, and set $X=X_1\cup X_2$. 
\item{$(a)$} $X$ is small $\Rightarrow$ $X_1$ and $X_2$ are both
small.
\item{$(b)$} $X$ is $2$-regular $\Leftrightarrow$ $X_1$ and $X_2$ are both $2$-regular.

\proof Part $(a)$ is easy.
To prove part $(b)$, write $L$ for the linear space
 $X_1\cap X_2$. The result follows from 
the exact sequence 
$$
0\rTo \I_X\rTo \I_{X_1}\oplus \I_{X_2}\rTo \I_{L}\rTo 0.
$$
(or see the more general statement in part $(c)$ of 
\ref{lin join equations}).\Box

\remark{} In the reduced case the converse to part $(a)$ follows from 
part $(b)$ and
\ref{main}.

\medskip

As with $2$-regularity, we can reduce questions about
linearly joined sequences to the connected case.

\proposition{direct sums} If $X\subset \PP^r$ is the
union of a linearly joined sequence of irreducible schemes,
then $X$ is the direct sum of its connected components.

\proof Let $X_1,\dots, X_n$ be the linearly joined
sequence of irreducible components of $X$.  We do induction on $n$,
the case $n=1$ being trivial.

By induction we may assume that $X'=\bigcup_{i=1}^{n-1}X_i$
is the direct sum of its connected components $Z'_1,\dots,Z'_s$.
Since $X_n\cap X'=\lspan(X_n)\cap \lspan(X')$
is a linear space, $X_n$ can meet at most one of the $Z'_i$.
If $X_n$ does not meet any $Z'_i$, and thus forms a new connected
component, then $\lspan (X_n)$ is disjoint from $\lspan(X')$, 
as required. Thus we may assume that $X_n$ meets a unique component, 
which we may as well call $Z'_s$, and the connected components of $X$ are
$$
Z_1=Z_1',\ \dots\ ,Z_{s-1}=Z_{s-1}',\ {\rm{and}}\  Z_s=Z_s'\cup X_n.
$$

If $X$ were not the direct sum of the $Z_i$, then there would
be a nontrivial dependence relation of the form
$\sum_1^{s-1} p_i+(p_s+q)=0$
where each $p_i$ is a vector of homogeneous coordinates of a point
in $\lspan(Z'_i)$, or the $0$ vector, and $q$ is a vector of homogeneous 
coordinates of a point in $\lspan(X_n)$. $q$ is not the $0$ vector since
$Z'_1,\dots,Z'_s$ are direct summands, so $q$ must represent a point in 
$\lspan(X')\cap \lspan(X_n)=X'\cap X_n\subset Z'_s$. Since
$X'$ is a direct sum of the $Z'_i$, it follows that 
$p_1=\cdots=p_{s-1}=p_s+q=0$, as required.\Box

Next we show that condition $(b)$ of \ref{small basic props} is 
exactly the difference between saying that the sequence of schemes
$X_1,\dots,X_n$ is linearly joined and saying that the (same) 
sequence of their spans is linearly joined.

\proposition{LJ for schemes and spans} Let  
$X_1,\dots,X_n\subset \PP^r$ be a sequence of closed subschemes and, 
for each $i$, let $L_i$ denote the linear span of $X_i$. The sequence
$X_1,\dots, X_n$ is linearly joined if and only if 
the sequence $L_1,\dots,L_n$ is linearly joined and each
pair $X_i, X_j$ is linearly joined, for all $i\neq j$.

\proof
First suppose that the sequence $X_1,\dots, X_n$ is linearly joined. It follows
at once that $L_1,\dots, L_n$ is linearly joined. By induction
on $n$ we may suppose that all pairs $X_i, X_j$ are linearly joined for
$i,j<n$, and it suffices to show that $L_j\cap L_n\subset X_j\cap X_n$
for each $j<n$. Since the sequence $X_1,\dots, X_n$ is linearly joined
we have
$$
\bigcup_{j<n}(X_j\cap X_n)=
(\bigcup_{j<n}X_j)\cap X_n=
\lspan(\bigcup_{j<n}X_j)\cap L_n.
$$
Since this is a linear space it must be 
contained in one of the $X_i\cap X_n$
for some $i<n$. In particular, $L_j\cap L_n\subset X_i\cap X_n$,
for all $j<n$, and thus we see that $L_j\cap L_n\subset X_n$. 
Since, by induction, the pairs $X_i,X_j$ are linearly joined 
we also have
$$
L_j\cap L_n= L_j\cap L_n\cap X_i\cap X_n\subset L_j\cap X_i
\subset X_j,
$$
completing the argument.

Conversely, suppose that the $X_i$ are pairwise linearly joined
and $L_1,\dots, L_n$ is a linearly joined sequence. By induction,
we may assume that $X_1,\dots, X_{n-1}$ is linearly joined sequence. 
But
$$
\lspan(\bigcup_{i<n}X_i)\cap L_n=
\lspan(\bigcup_{i<n}L_i)\cap L_n=
(\bigcup_{i<n}L_i)\cap L_n= L_j\cap L_n
$$
for some $j<n$. Since the pair $X_j,X_n$ is linearly joined we also 
have $L_j\cap L_n=X_j\cap X_n$, so
$\lspan(\bigcup_{i<n}X_i)\cap L_n=(\bigcup_{i<n}X_i)\cap X_n$
as required.\Box

\section{basic properties} Proof of \ref{main}

We already have the tools to dispose of two implications in \ref{main}
easily. \ref{homol char of min deg}, together with
\ref{linearly joined equivalences} gives the implication $(c)\Rightarrow (b)$.  
On the other hand, \ref{2-regular gives small} (proved in Eisenbud-Green-Hulek-Popescu [2004]) 
includes the  implication $(b) \Rightarrow (a)$ as a special case.

The last implication, $(a)\Rightarrow (c)$,
will occupy us for the  rest of this section. Here is an outline:
We first prove $(a)\Rightarrow (c)$ in the case where each
$X_i$ is a linear space, then we use this case to 
prove the implication $(b)\Rightarrow (c)$ in general. 
Finally we will use the implication  
$(b)\Rightarrow (c)$ to prove $(a)\Rightarrow (c)$.

\medskip
\noindent{\it Proof of $(a)\Rightarrow (c)$ for unions of planes}:
\smallskip

To finish the proof
of \ref{main} in the case where $X$ is a union of planes,
we will project
from a general point on a carefully chosen
component of $X$, using the following results.

\proposition{matroid} Let $V_i\subset V$, $i=1,\ldots,m$, 
be distinct linear subspaces, not contained in
one another,  such that any $m-1$ of them span
the ambient space $V$. Then there exist vectors $u_i\in V_i\setminus
\bigcup_{j\neq i}V_j$ such that the collection $u_1,\ldots, u_m$ is 
linearly dependent.

\proof Since the subspaces $V_i$ are distinct and not contained in
one another $Z_i=(\cup_{j\neq i}V_j)\cap V_i$ is
a union of finitely many proper linear subspaces of $V_i$,
for all $i=1,\ldots, m$. In particular $Z_i\neq V_i$.
Let $W$ denote the kernel of
the natural summation map
$$
\varphi:V_1\oplus V_2\oplus\cdots\oplus V_m\rTo^{(v_1,\ldots,v_m)
\mapsto\sum_{i=1}^m v_i} V.
$$
Observe that if $\pi_i:V_1\oplus V_2\oplus\cdots\oplus V_m\to V_i$
denotes the projection on the $i$-th factor, then $\pi_i(W)=V_i$ by
our hypothesis.  It follows that $\pi_i^{-1}(Z_i)\neq W$
for all $i$. Since the ground field is infinite, there exists a vector
$w=(v_1,\ldots,v_m)\in W$ such that $\pi_i(w)\not\in Z_i$ for all
$i=1,\ldots,m$. In other words, $v_i\in V_i\setminus(\cup_{j\neq
i}V_j)$, $i=1,\ldots,m$ and $v_1+\cdots+v_m=0$. \Box

\proposition{simple secants} Let $X=\bigcup_{i=1}^t X_i\subset\PP^r$ 
be a small union of closed subschemes. If $\Lambda_i\subset X_i$ is a linear subspace
for each $i$ such that $\Lambda_i$ does not meet any $X_j$ for
$j\neq i$, then the $\Lambda_i$ are linearly independent and their
linear span meets each $X_j$ precisely in $\Lambda_j$.

\proof Let $\Lambda$ be the span of the $\Lambda_i$.
If the $\Lambda_i$ were dependent, then there would be
a set of points $p_i\in \Lambda_i$ that were dependent. Similarly
if $\Lambda\cap X_j\neq \Lambda_j$ then there would be a set of 
points $p_i\in \Lambda_i$, such that the span of the $p_i$ met
$X_j$ outside $\Lambda_j$. Consequently, to prove either
statement, we may assume that $\Lambda_i=\{p_i\}$.

Since $p_i\notin X_j$ for $j\neq i$, there is a component $X_i'$ of 
$\Lambda\cap X$ containing $p_i$ but none of the $p_j$. 
By \ref{small char}, the sum of the degrees of the $X_i'$
is at most $\codim(X\cap \Lambda,\Lambda)+1\le\dim(\Lambda)+1$. Since
$\Lambda$ is spanned by the points $p_i\in X_i'$, the dimension
of $\Lambda$ is at most one less than the number of components
$X_i'$. By \ref{small char} we conclude that each $X_i'=\{p_i\}$, 
the $p_i$ are linearly independent, and there
are no other points in $X\cap \Lambda$.\Box

\corollary{small matroid}
If $X=\bigcup_i X_i\subset\PP^r$ is a small union of linear subspaces, then there
 exists an $i$ such that $X_i$ is not in the linear span of the
 $\bigcup_{j\neq i}X_j$.

\proof If each component $X_i$ were contained in the 
span of the others, then passing to affine cones we would
have a system of subspaces satisfying the hypothesis of
\ref{matroid}. The   set of points corresponding to the 
vectors in the conclusion of \ref{matroid} would 
violate \ref{simple secants}.\Box

\proposition{small projects small} Let $X=\bigcup_{i=1}^nX_i\subset \PP^r$ be
a small union of linear subspaces.  If $p\in X_1$ is a point that 
is not on the span of $\bigcup_{i=2}^n X_i$, then 
the image of $X$ under linear projection from $p$, 
$\pi_p:\PP^r\rDotsto \PP^{r-1}$, is small.

\proof By hypothesis there is a hyperplane $H\subset \PP^r$ that
contains $\bigcup_{i=2}^n X_i$ but does not contain $p$. By \ref{small basics}
the union of linear spaces $X'=(H\cap X)_\red$ is small,
and the projection $\pi_p$ induces a linear isomorphism
$X'\subset H\rTo \pi_p(X)\subset \PP^{r-1}$, proving that
$\pi_p(X)\subset \PP^{r-1}$ is also small.\Box

\proposition{unproject linearly joined}
Let $X=\bigcup_{i=1}^n X_i\subset\PP^r$ be a small union of linear spaces,
and let  $p\in X_1$ be a point not in the linear
span of the union of  $X_2,\dots, X_n$.
If, for some $i$, the schemes
$\pi_p(X_i)$ and $\pi_p(\bigcup_{j\neq i}X_j)$ are linearly joined,
then $X_i$ and $\bigcup_{j\neq i}X_j$ are  linearly joined.

\proof Let $Z=\bigcup_{j\neq i}X_i$. We must show that if 
$q\in X_i\cap \lspan(Z)$ then $q\in Z$. Note
that $p$ cannot be in both $X_i$ and $\lspan(Z)$, so $p\neq q$,
and the projection $\pi_p(q)$ is defined.

First suppose $i=1$. We have
$\pi_p(q)\in\pi_p(X_1\cap \lspan(Z))=
\pi_p (X_1)\cap \pi_p(\lspan(Z))$
because $p\in X_1$. Moreover, $\pi_p(\lspan (Z))=\lspan (\pi_p(Z))$,
and thus, because  $\pi_p(X_1)$ and $\pi_p(Z)$
are linearly joined,
$$
\pi_p(q) \in \pi_p(X_1)\cap \lspan (\pi_p(Z))= \pi_p(X_1)\cap \pi_p(Z).
$$
In particular, $\pi_p(q)\in \pi_p(Z)$.
As $q\in \lspan(Z)$, and 
$\pi_p$ is an isomorphism on $\lspan(Z)$, we get
$q\in Z\cap X_1$ as required.

Next suppose $i\neq 1$. Because $p\in \lspan(Z)$, we may argue
as before and obtain
$$
\pi_p(q)\in\pi_p( X_i )\cap\pi_p(\lspan(Z))=
\pi_p(X_i)\cap\pi_p(Z).$$
Thus the line $\lspan(p,q)$ meets $Z$ in a point $q'$ and 
meets $X_i$ in a point $q''$.
If $q=q'$, then $q\in Z$, and we are done. Thus we may as well assume
that $q\neq q'$, which implies that $p\in \lspan(q,q')$.
By hypothesis $p\notin \lspan(\bigcup_{i\neq 1}X_i)$, 
so at least one of the points $q, q'$ must be in $X_1$. 
Since  $p\in X_1$, this means that both $q,q'\in X_1\subset Z$; 
in particular $q\in Z$ as required.\Box

\noindent{\it Conclusion of the Proof of $(a)\Rightarrow (c)$
for a union of linear spaces.\/}
Again, let $X=\bigcup_{i=1}^n X_i\subset\PP^r$ be a small
union of linear spaces. We do induction on 
the dimension $r$ of the ambient projective space. 
After renumbering the components we may assume by
\ref{small matroid} that there exists a point $p\in X_1$ that
is not in the linear span of $\bigcup_{i\geq 2}X_i$. 
By \ref{small projects small}, the
projection $\pi_p(X)\subset\PP^{r-1}$ is also small. 

By 
induction on the dimension of the ambient space, 
there is a permutation of $\{1,\ldots,n\}$, 
such that $\pi_p(X_{\sigma(t+1)})$ is linearly joined to  
$\bigcup_{i=1}^t\pi_p(X_{\sigma(i)})$
for $t=1,\dots,n-1$. If $1\notin\{\sigma(1),\dots, \sigma(t+1)\}$,
then $\pi_p$ is an isomorphism on the linear span of 
$\bigcup_{i=1}^{t+1}X_{\sigma(i)}$, so $X_{t+1}$ is linearly
joined to  $\bigcup_{i=1}^{t}X_{\sigma(i)}$. If on the
other hand $1\in\{\sigma(1),\dots, \sigma(t+1)\}$,
then \ref{unproject linearly joined} yields the same conclusion.
 \Box

\medskip
\noindent{\it Proof of $(b)\Rightarrow (c)$ in the general case.\/}
If $X$ is  $2$-regular, then by \ref{small basic props} all irreducible
 components $X_i$ of $X$ are varieties of minimal degree
in their spans, and they are pairwise linearly joined.
To prove $(c)$ it suffices, by \ref{LJ for schemes and spans}, to show
that the union of the linear spans of the $X_i$ can be
arranged in a linearly joined sequence. The next result shows
that this union is $2$-regular, and thus
reduces the proof to the case we have already treated.
The proof will use an induction, and for this we need to
know that taking hyperplane sections commutes with taking
spans in certain cases.

\lemma{spans and sections} Let $Y\subset \PP^s$ be a scheme,
and $H$ the hyperplane defined by a linear form $h$.  Suppose $h$ is
sufficiently general that it is locally a nonzerodivisor on $\O_Y$. If
$h^0(\O_Y)=1$, then $\lspan(H\cap Y) = H\cap\lspan(Y)$.

\proof
The diagram 
$$\diagram[small,midshaft]
0&\rTo& H^0(\O_{\PP^s})&\rTo^h& H^0(\O_{\PP^s}(1))&\rTo^\rho &H^0(\O_{H}(1))&\rTo 0\cr
  &&\dTo&&\dTo&&\dTo\cr
0&\rTo& H^0(\O_Y)&\rTo^h& H^0(\O_Y(1))&\rTo &H^0(\O_{Y\cap H}(1))&\rTo&\cdots
\enddiagram
$$
has exact rows, and the left-hand vertical map is surjective by the 
connectedness of $Y$. By the snake lemma, the restriction $\rho$
induces the isomorphism $H^0(\I_{Y}(1))\cong H^0(\I_{Y\cap H,H}(1))$.  
It follows that $\lspan(H\cap Y) = H\cap\lspan(Y)$.\Box

\theorem{oberwolfach 2-regular} Let $X=\bigcup_iX_i\subset\PP^r$ 
be a closed subscheme with irreducible components $X_i$ that
are Cohen-Macaulay and $2$-regular. 
Let $L_i=\lspan(X_i)$ and set $Y=\bigcup_i L_i$.
If $X_i\cap X_j =L_i\cap L_j$ for each $i\ne j$, then $X$ is
$2$-regular if and only if
$Y$ is $2$-regular.

\proof
\noindent We will show that conditions $(a)$ and $(b)$
of \ref{2-reg criterion} hold for $X$ if and only if they
hold for $Y$. 

\smallskip\noindent $(a)$:
Suppose that $X$ is $2$-regular. By \ref{direct sums 2-reg}, $X$ is the direct sum of
its connected components, which are also $2$-regular.
In particular, if $X$ has any zero-dimensional components,
they form a direct summand and we may drop them. Thus we
may assume that every component of $X$ has dimension $\geq 1$.

Let $H\subset \PP^r$ be a general hyperplane, so that that $H\cap X$ is $2$-regular.
Since $X_i$ has dimension $\geq 1$ and is Cohen-Macaulay,
$h^0(\O_{X_i})=1$. By \ref{spans and sections},  $H\cap L_i$ 
is the span of $H\cap X_i$. 

By \ref{direct sums 2-reg} again, $H\cap X$ is the direct sum of
its connected components, which are also $2$-regular,
and it suffices to show that the same is true of $H\cap Y$.
The connected components of $H\cap Y$ have the form
$\bigcup_{i\in I}H\cap L_i$ where
$I$ is a minimal set of indices $i$ such that
$$
\dim\bigl[\bigcup_{i\in I}L_i \cap \bigcup_{j\notin I} L_j\bigr]\le 0.
$$
By hypothesis $L_i\cap L_j=X_i\cap X_j$, for all $i\neq j$,
so 
$\bigcup_{i\in I}H\cap X_i$ is a union of connected components
of $H\cap X$ (extra components appear when there are $1$-dimensional
components of $X$). Since such a union is a direct summand of $H\cap X$,
it follows that $H\cap Y$ is a direct sum of its connected components.

Each $1$-dimensional irreducible component $X_i$ 
of $X$ contributes a zero-dimensional scheme spanning a 
single connected component $H\cap L_i$ of $H\cap Y$ that is a single
linear space, and is thus $2$-regular. The union of the higher-dimensional
irreducible  components of $X$ contributes the union of a subset of the 
connected components of $H\cap X$, which is thus $2$-regular.
But for an irreducible component $X_i$ of dimension at least $2$,
$H\cap X_i$ is again irreducible and spans $H\cap L_i$, by \ref{spans and sections}.
Thus we may use induction, and deduce that $H\cap Y$ is
$2$-regular as required.

The case where $H\cap Y$ is $2$-regular is similar. Since it will
not be required for the proof of \ref{main}, and is a consequence
of that result in the case where $X$ is reduced, we omit
the details.

\smallskip\noindent $(b)$:
Consider the diagram
$$\diagram
0&\rTo& H^0(\O_Y(1))&\rTo& \bigoplus_i H^0(\O_{L_i}(1))
&\rTo& \bigoplus_{i<j} H^0(\O_{L_i\cap L_j}(1))\cr
&&\dTo^{\rho_{Y/X}}&&\dTo^{\oplus\rho_{L_i/X_i}}
&&\dEq_{\oplus\rho_{L_i\cap L_j/X_i\cap X_j}}\cr
0&\rTo& H^0(\O_X(1))&\rTo& \bigoplus_i H^0(\O_{X_i}(1))
&\rTo& \bigoplus_{i<j} H^0(\O_{X_i\cap X_j}(1))\cr
\enddiagram
$$
where the vertical maps are restrictions. The maps
$\rho_{L_i/X_i}$ are isomorphisms because the 
$X_i$'s are linearly normal and each span $L_i$.
Since $X_i\cap X_j=L_i\cap L_j$ the right hand side vertical
map is an equality. Thus the map $\rho_{Y/X}$ is also an isomorphism,
so $X$ and $Y$ are either both linearly normal in $\PP^r$ or both not. 
\Box

\medskip
\noindent{\it Proof of \ref{main} continued.\/}
We complete the proof of \ref{main}
by proving $(a)\Rightarrow (b)$, using induction
on $\dim X$. A general hyperplane section of $X$ is reduced by 
Bertini's Theorem, and by \ref{small basics}, it is small.
By induction, it is $2$-regular. Thus
by \ref{2-reg criterion}, it is enough to show that
$X$ is linearly normal.

If $X$ were not linearly normal, we could write $X$ as
a linear projection of a linearly normal 
variety $Y$ in some larger
projective space $\PP^s=\PP(H^0(\O_X(1))$ from 
a linear space $M\subset \PP^s$ that is disjoint from $Y$
in such a way that the projection is an isomorphism on $Y$,
but $M$ is contained in the linear span of $Y$.
We first want to show that $Y\cup M$ is small. For this
we will need the following basic fact about geometric degree:

\proposition{iso projection and plane section} Suppose that
$Y\subset \PP^s$ is a scheme, and $M\subset \PP^s$ is 
a linear space disjoint from $Y$. Suppose that the linear projection
from $M$, denoted $\pi_M:  \PP^s\rDotsto \PP^r$, induces an 
isomorphism $Y\to \pi_M(Y)$. If $L\subset \PP^s$ is a linear space containing $M$,
then $\pi_M(L)\cap \pi_M(Y) = \pi_M(L\cap Y) \cong L\cap Y$,
and the geometric degree of $L\cap Y\subset \PP^s$ is the same as
the geometric degree of $\pi_M(L)\cap \pi_M(Y)$.

\proof First note that because $M\subset L$ the map
$\pi_M(Y\cap L)\to \pi_M(Y)\cap \pi_M(L)$ is 
set-theoretically surjective, so it is
enough to check the result locally at a point $q\in Y\cap L$.

Suppose that $W$ is the space of linear forms vanishing on
$M$, and that $V\subset W$ is the space of linear forms vanishing
on $L$. We may identify $W$ with $H^0(\O_{\PP^r}(1))$, and under this 
identification, $V$ is identified with the space of linear forms on 
$\PP^r$ that vanishes on $\pi_M(L)$. With this understanding, the first
claim in the  statement of the theorem may be written, locally at  $q$, as 
$$
\O_{Y\cap L, q}\cong \O_{Y,q}/V\O_{Y,q} 
\cong \O_{\pi_M(Y),\pi_M(q)}/V\O_{\pi_M(Y),\pi_M(q)},
$$
which holds because $\pi_M$ induces an isomorphism 
$\O_{Y,q}\cong \O_{\pi_M(Y),\pi_M(q)}$.

The statement about geometric degrees follows because all the 
isomorphisms preserve linear forms. \Box

\proposition{unproject some small} Suppose that $Y\subset \PP^s$
is an algebraic set, and that $M\subset \PP^s$ is a linear
space, disjoint from $Y$, such that the linear projection
$\pi_M$ induces an isomorphism $Y\to \pi_M(Y)$. If $\pi_M(Y)$
is small, then so is $Y\cup M$.

\proof Suppose that $\Lambda\subset \PP^s$ is a plane that meets
$X=Y\cup M$ in a finite set. We must show that 
$\deg(X\cap \Lambda)\leq 1+\dim(\Lambda)$. Let $L=\lspan(\Lambda\cup M)$. 
As $Y\cap\Lambda$ is a plane section of $Y\cap L$, the degree of $Y\cap\Lambda$ 
is bounded by the geometric degree of $Y\cap L$. 
\ref{iso projection and plane section} implies that this geometric degree 
is the same as the geometric degree of 
$\pi_M(Y\cap L)=\pi_M(Y)\cap\pi_M(L)=\pi_M(Y)\cap\pi_M(\Lambda)$. 
Since $\pi_M(Y)$ is small by hypothesis, this is bounded by 
$\dim(\pi_M(\Lambda))+1$. Combining this inequality with
the obvious $\deg(M\cap\Lambda)+\dim(\pi_M(\Lambda))=\dim(\Lambda)$,
we deduce that $\deg(X\cap\Lambda)\leq 1+\dim(\Lambda)$ as required.\Box

\smallskip
\noindent{\it Proof of \ref{main} continued.\/}
{}From the fact that $Y\cup M$ is small, it follows by
\ref{subtracting disjoint set} that  $Y$ is small. As $Y$ is linearly normal
by construction, we see from  induction on the dimension, \ref{small basics} $(a)$, 
and \ref{2-reg criterion} that $Y$ is $2$-regular. We want to
show next that $Y$ actually coincides with $X$ (i.e. $M$ is
empty), and as a first step we show the following result
which interesting on its own. 

\proposition{2-regular secants}
Let $Y\subset \PP^s$ be a $2$-regular algebraic set.
If $p\in \PP^s$ is a point in the span of $Y$,
 then there is a plane $L$ containing $p$
such that $L\cap Y$ is finite and $L\cap Y$ spans $L$.

\proof  Using the implication
$(b)\Rightarrow (c)$ of \ref{main}, we see that $Y$ is the union
of a linearly joined sequence of irreducible varieties of
minimal degree $Y_1,\dots, Y_n$. By \ref{small basic props} $(a)$
and  \ref{2-reg of min deg} $(c)$, the homogeneous 
coordinate rings of the irreducible components of $Y$
are Cohen-Macaulay. We will prove, by induction on
the dimension of $Y$, a more general result: if $Y\subset \PP^s$
is the union of a sequence of linearly joined irreducible schemes
$Y_1,\dots Y_n$ such that each homogeneous coordinate
ring $S_{Y_i}$ is Cohen-Macaulay, and no $(Y_i)_\red$
is contained in another, then for any point
$p\in \lspan(Y)$, there exists a plane $L$ containing $p$
such that  $L\cap Y$ is finite and $L\cap Y$ spans $L$.

By \ref{direct sums} the problem reduces to the case where 
 $Y$ is connected. If $Y$ is supported
at a point the assertion is obvious. Otherwise, each
component of $Y$ must have dimension $\geq 1$.
Under this circumstance, we show that $h^0(\O_Y)=1$. 
Since $S_{Y_i}$ has depth at least $2$ we have
$
\bigoplus_m H^0(\O_{Y_i}(m))\cong S_{Y_i},
$
so $h^0(\O_{Y_i})=1$ for every $i$. In general,
if $A, B\subset \PP^s$ are schemes with $h^0(\O_A)=h^0(\O_B)=1$,
and if $A\cap B\neq \emptyset$, then taking the cohomology of
the short exact sequence
$$
0\rTo \O_{A\cup B}\rTo \O_A\oplus \O_B\rTo \O_{A\cap B}\rTo 0
$$
gives that $h^0(\O_{A\cup B})=1$ also. In particular,
$h^0(\O_Y)=1$ as claimed.

Let $H$ be a general hyperplane  containing $p$,
defined by a linear form $x\in H^0(\O_{\PP^s}(1))$.
Since $H$ cuts properly every component $Y_i$ of $Y$,  
$S_{H\cap Y_i} = S_{Y_i}/(x)$ is again
Cohen-Macaulay. It follows from the definitions
that the sequence $H\cap Y_1, \dots ,H\cap Y_n$ is linearly joined.
Since $Y$ is connected and reduces we see that $h^0(\O_Y)=1$. By
\ref{spans and sections}
$H$ is spanned by $H\cap Y$, and we are done.\Box

\corollary {not small}
If $Y$ is a $2$-regular algebraic set, and $Y'$ is a scheme disjoint
from $Y$ that meets the span of $Y$, then $Y\cup Y'$ is not small.

\proof By \ref{2-regular secants} there is
a linear space $L$ that contains a point of $Y'$ and meets $Y$ in a
finite scheme that spans $L\cap Y$.  Whatever the dimension of $L\cap
Y'$ this violates the conclusion of \ref{small char}, 
so $Y\cup Y'$ is not small.\Box

\medskip
\noindent{\it Conclusion of the proof of \ref{main}.\/}
We may apply \ref{not small} to the case of the small
schemes $Y\cup M$ from \ref{unproject some small} and
$Y$ to conclude that $M$ is empty! That is, the original
small scheme $X$ was linearly normal, and thus $2$-regular.
This finishes the proof of $(a)\Rightarrow (b)$,
and with it the proof of \ref{main}.\Box

\remark{} The conclusion of \ref{2-regular secants}
seems to be true for a wide class of schemes $Y$. We are grateful to
Harm Derksen for pointing out that it does not hold for arbitrary
schemes, however: Let $M_1\subset M_2\subset \PP^4$ be a line
contained in a $2$-plane in $\PP^4$, and let $p_1, p_2$ be general
points of $M_2$.  Let $Y=M_1\cup P_1\cup P_2$, where $P_i$ is a double
point with $(P_i)_\red=p_i$ and general tangent vector. It is
immediate that $Y$ spans $\PP^4$. However, if $q$ is a general point
of $\PP^4$, then $q\notin \lspan(M_1\cup P_1)\cup
\lspan(M_1\cup P_2)\cup\lspan(P_1\cup P_2)$
and so any linear space $L$ containing 
$q$ and spanned by $L\cap Y$ must contain at least $p_1, p_2,$ and a 
point of $M_1$ which is not collinear with $p_1, p_2$. It follows easily 
that such an $L$ is all of $\PP^4$, and thus
$L\cap Y$ is not finite. 

On the other hand, it seems that the conclusion of 
\ref{2-regular secants} might hold for all reduced schemes.  
Kristian Ranestad has proven this in characteristic $0$.  
We are grateful to him for allowing us to include it here.

\proposition{ranestad} Let $X\subset\PP^s$ be a reduced
scheme over an
algebraically closed field of characteristic zero.
If $p$ is a point in the linear span of $X$, then  
there is a linear subspace $L\subset\PP^s$ containing $p$ such that
$L\cap X$ is finite and contains a set of distinct points that spans
$L$.

\proof  We may as well assume that $X$ spans $\PP^s$. 
If $p\in X$ or $X$ is finite there is nothing to prove. So let $X_0$
be a positive dimensional irreducible component of $X$ and let
$L_0$ be the linear span of $X_0$. Since $L_0$ has dimension $r>0$, 
we may choose distinct points $x_0,\ldots,x_r\in X_0$, and 
$x_{r+1},\ldots,x_s\in X\setminus X_0$, such that $x_0,\ldots,x_s$ 
span $\PP^s$. 

If $L\subset \PP^s$ is a linear space of dimension $<s$ such that $L$
is spanned by some set of distinct points of $X$ and $p\in L$, then we
may replace $X$ by $(L\cap X)_\red$ and we are done by induction on
$s$.  Thus for example we are done if $p\in L_1=\lspan\{x_{r+1},\ldots,x_s\}$.  
Therefore we may assume that $p\notin L_1$, so the plane $L_2=\lspan(p,L_1)$ 
meets $L_0$ in a unique point $q$.  If $q\in X$, then since
$L_2=\lspan\{q,x_{r+1},\ldots,x_s\}$ and $r\geq 1$, we are again done
by induction.

Thus we may suppose $q\not\in X$. It follows that
the system  of hyperplanes containing $L_2$ has no basepoints on $X_0$. 
Since the the base field is of characteristic zero it follows that a general 
hyperplane $H$ containing $L_2$ meets $X_0$
in a reduced set. Since $X_0$ is reduced and
irreducible, \ref{spans and sections}
shows that $H\cap X_0$ spans $H\cap L_0$.
Therefore we may find distinct points $y_1,\ldots,y_r\in H\cap X_0$ 
such that $y_1,\ldots,y_r,x_{r+1},\ldots,x_s$ span $H$.
Since $p\in H$ again we are done by induction on $s$, and the proposition
follows.\Box

\section{chordal} Small subspace arrangements

In this section we will study the combinatorics of the condition for
union of linear subspaces to be small.  Recall that a subgraph $F$ of
a graph $G$ is called a {\it forest\/} if $F$ is a disjoint union of
{\it trees\/} (acyclic connected graphs). A {\it leaf\/} of $F$ is a
vertex of $F$ connected to at most one other vertex of $F$.  A forest
$F\subset G$ {\it spans\/} $G$ if $F$ contains all the vertices of $G$
and any two vertices connected by a path in $G$ are connected by a
path in $F$.  We say that an ordering of the vertices of $G$ is {\it
compatible with a spanning forest $F\subset G$\/} if the smallest
vertex in every connected component is a leaf, and the ordering
restricts to the natural ordering on the vertices of any path starting
from that leaf.

By a {\it subspace arrangement\/} we mean a finite union of
incomparable linear subspaces in a projective space, say
$Y=\bigcup_{i=1}^nL_i\subset \PP^r$. We generally do not distinguish
between the set of subspaces and their union. To a subspace
arrangement $Y$ we associate the weighted graph $G_Y$, whose vertices
are the subspaces $L_i$ of $Y$, and whose edges join the pairs of
subspaces with non-empty intersection.  We define the {\it weight\/}
of the edge $(L_i, L_j)$ to be $1+\dim(L_i\cap L_j)$, and the weight
of a subgraph is the sum of the weights of its edges. We will be
interested in the spanning forests of maximal weight in $G_Y$. To
simplify the notation, we give the vertex $L_i$ weight $1+\dim(L_i)$
and, for any graph $G$ with edge and vertex weights 
we define the {\it weighted Euler characteristic\/} $\chi_w(G)$ 
to be the sum of the weights of the vertices minus the sum of 
the weights of the edges. Thus a spanning forest $F$ for $G_Y$ 
has maximal weight if and only if it has minimal $\chi_w(F)$.

The main result of this section says that smallness 
and linearly joined sequences of components of $Y$ are
features of certain spanning forests in $G_Y$.
\goodbreak

\theorem{weights of forests} Let $Y\subset \PP^r$ be
a subspace arrangement.
\item{$(a)$}   $Y$ is small if and only if
the weighted graph $G_Y$ has a spanning forest $F$ with
$\chi_w(F)= 1+\dim (\lspan(Y)),$ the smallest possible value.
\item{$(b)$} If $Y$ is small, then an ordering of the components of $Y$ 
makes them into a linearly joined sequence if and only if it is 
compatible with some spanning forest satisfying the equality above.

We begin with an elementary result that shows the above given value of
$\chi_w$ is the smallest possible and explains the connection with
linearly joined sequences.

\lemma{dimension of vector spaces} 
 Let $Y=\bigcup_{i=1}^n L_i\subset \PP^r$ be a subspace arrangement, and let
$F$ be a spanning forest of $G_Y$. 
\item{$(a)$} Suppose that $L_n$ is a leaf of $F$. If
$Y'$ and $F'$ are obtained from $Y$ and $F$ by removing $L_n$,
then 
$$
\chi_w(F)\geq \chi_w(F')+\bigl[\dim(\lspan(Y))-\dim(\lspan(Y'))\bigr],
$$
with equality if and only if $L_n$ is linearly joined with
$Y'$ and either $L_n\cap Y'$ is empty or $L_n\cap Y' = L_n\cap L_j$, where
$L_n$ is connected to $L_j$ in $F$.
\item{$(b)$} $\chi_w(F)\geq 1+\dim(\lspan(Y))$.

\noindent {\it Proof of \ref{dimension of vector spaces}.\/}
 Part $(a)$ is elementary, and part $(b)$ follows from
part $(a)$ by induction on the number of components.\Box

\noindent {\it Proof of \ref{weights of forests}.\/}
We prove parts $(a)$ and $(b)$ together.
First suppose that $F$ is a  spanning forest with
$\chi_w(F)= 1+\dim (\lspan(Y))$, and
$L_1,\dots, L_n$ is a compatible ordering of the 
subspaces in $Y$. 

It follows that $L_n$ is a leaf
of $F$. Let $Y'$ and $F'$ be obtained by deleting
$L_n$ from $Y$ and $F$ respectively. By 
\ref{dimension of vector spaces} we have
$$\eqalign{
1+\dim(\lspan(Y'))
&\leq
\chi_w(F')\cr
&\leq \chi_w(F)-\bigl[\dim(\lspan(Y))-\dim(\lspan(Y'))\bigr]\cr
&=1+\dim(\lspan(Y'))
}
$$
so all the equalities hold. Thus $L_n$ is linearly
joined to $Y'$ and, by induction on the number of components,
$L_1,\dots, L_n$ is a linearly joined sequence
and $Y'$ is small. Using \ref{main} it 
follows from \ref{linearly joined equivalences} 
that $Y$ is small.

Now suppose that $Y$ is small. By \ref{main} $(c)$ we may order the
components of $Y$ to form a linearly joined sequence $L_1,\dots,
L_n$. It suffices to find a spanning forest with which this ordering
is compatible. We do induction on the number $n$ of components of $Y$,
the case $n=1$ being trivial.

Since the sequence $L_1,\dots, L_{n-1}$ is linearly joined,
$Y'=\bigcup_{i=1}^{n-1} L_i$ is also small. By
induction, $G_{Y'}$ contains a  
spanning forest $F'$ with $\chi_w(F')=1+\dim(\lspan(Y'))$
such that the ordering $L_1,\dots,L_{n-1}$ is compatible
with $F'$.
Since $Y'\cap L_n=\lspan(Y')\cap L_n$ 
is a linear space, and the ground field is infinite,
$Y'\cap L_n$ is either empty
or it is equal to $L_j\cap L_n$ for some $j<n$. In 
the first case $L_n$ is a connected component of $G_Y$,
so adjoining it to $F'$ we get a spanning forest $F$ of $G_Y$ 
In the second case, we may adjoin $L_n$ to $F'$
and connect it to $L_j$, obtaining a new spanning forest.
In either case we get equality in the formula of
\ref{dimension of vector spaces} $(a)$, and the 
given order is compatible with the forest $F$.\Box

\ref{weights of forests} makes it interesting to
understand better which forests have minimal
weighted Euler characteristic. If $Y=\bigcup L_i$ is a 
subspace arrangement, we say that a forest $F\subset G_Y$ satisfies
the {\it clique-intersection property\/}
if whenever $L_1,\dots, L_j$ form a path in
$F$ we have
$$
L_1\cap L_j=L_1\cap L_2\cap\cdots\cap L_j.
$$

\proposition{tarjan}  Let   
$Y=\bigcup_{i=1}^n L_i\subset\PP^r$ be a small subspace arrangement
with intersection graph $G_Y$.  A spanning forest $F\subset G_Y$ has
$\chi_w(F)=1+\dim(\lspan(Y))$ (or equivalently, maximal weight) if and
only if $F$ satisfies the clique-intersection property.

\proof 
Prim's algorithm [1957] (see also Graham-Hell [1985]) shows
that a spanning forest $F$ in a weighted connected 
graph $G$ has minimal $\chi_w(F)$ 
if and only if, for each edge $(x,y)$ of $G$,
the path in $T$ joining $x$ to $y$
consists of edges (each) of weight $\geq$ the weight of 
$(x,y)$ (see also Tarjan [1983, pp.~71--72].)
In particular, a spanning forest 
$T\subset G_Y$ satisfying the  clique-intersection property 
must have minimal weighted Euler characteristic.
Since $Y$ is small, \ref{dimension of vector spaces}
shows this is $1+\dim(\lspan(Y))$.

On the other hand, suppose that $F$ is any spanning forest with
$\chi_w(F)=1+\dim(\lspan(Y))$
and $L_1,\dots, L_j$ are the spaces
along a path in $F$.  By \ref{weights of forests}
the spaces $L_1,\dots, L_j$ form a part of
a linearly joined  sequence 
$L_{-t},\dots, L_0, L_1,\dots, L_j, L_{j+1},\dots, L_s$
involving all
the components of $Y$. It follows that
all of the $L_i\cap L_j$, for $i<j$ are contained in 
one $L_h\cap L_j$, with $h<j$.
By Prim's algorithm, $\dim(L_h\cap L_j)\leq \dim(L_{j-1}\cap L_{j})$
so all the $L_i\cap L_j$ are contained in $L_{j-1}\cap L_j$.
Thus by induction on $j$ we get
$$\eqalign{
L_1\cap L_j&= L_1\cap L_{j-1}\cap L_j\cr
&=
L_1\cap L_2\cap\ldots \cap L_{j-1}\cap L_j
}
$$
as required. \Box

Another way of finding an ordering of the components of a small
subspace arrangement $Y$ as a linearly joined sequence is the
following: select vertices of $G_Y$ inductively by choosing, at each
step, an unselected vertex $i_k$ which has the maximal number of
adjacent vertices among the vertices already selected. For a proof see
Tarjan-Yannakakis [1984].

The union of all spanning forests of maximal weight of $G_Y$ 
is the subgraph $H_Y\subset G_Y$ with the same vertices as the intersection 
graph $G_Y$, but where $L_i$ and $L_j$ are joined by an edge only
when $L_i\cap L_j$ disconnects $Y$. Indeed, by \ref{tarjan} a
maximal weight spanning forest $F$ of $G_Y$ satisfies the 
clique-intersection property which easily implies that $F$
is actually a spanning forest in $H_Y$. Conversely, if 
$L_i\cap L_j\ne\emptyset$ disconnects $Y$ and 
$L_i=L_{k_0},L_{k_1},\ldots, L_{k_r}=L_j$ is the path joining $L_i$ and $L_j$ 
in a maximal weight spanning forest $F$ of $G_Y$, then necessarily
$L_{k_m}\cap L_{k_{m+1}}= L_i\cap L_j$, for some $m<r$ and the
forest $F'$ obtained by replacing in $F$ the edge $(L_{k_m},L_{k_{m+1}})$
by $(L_i,L_j)$ is also of maximal weight.

\section{equations} Equations and syzygies of reduced 2-regular schemes

Next we show how to find generators for the ideal of a reduced
$2$-regular projective scheme using property $(c)$ of \ref{main}.
Consider a closed subscheme $X=X'\cup X''\subset\PP^r$, and set $L'=\lspan(X'),\
L''=\lspan(X'')$.  In general it is difficult to find generators for
the intersection of two ideals, but if $X'$ and $X''$ are linearly
joined (that is, $X'\cap X''=L'\cap L''$) then we can give minimal
generators and a (non-minimal) free resolution of $I_X=I_{X'\cup
X''}=I_{X'}\cap I_{X''}$ explicitly from minimal generators and free
resolutions for $I_{X'/L'}$ and $I_{X''/L''}$.  Our result extends
results of Barile and Morales [2000, 2003].

For simplicity we will suppose throughout that $L'\cup L''$ spans the
whole ambient space $\PP^r$, and leave the reader the easy task to
adapt \ref{lin join equations} below to the degenerate case.  We write
$\mu(I)$ for the minimal number of generators of a homogeneous ideal
$I$, and $\reg(I)$ for its regularity.

\theorem{lin join equations} 
Let $X=X'\cup X''\subset \PP^r$ be a nondegenerate closed scheme that
is the union of two subschemes $X'$ and $X''$ with linear spans $L'$
and $L''$, respectively. Suppose that $X'$ and $X''$ are linearly
joined along $L=L'\cap L''$.
\item{$(a)$}
$$
I_X = 
\widetilde {I_{X',L'}}+
\widetilde{I_{X'',L''}}+
I_{L'}\cdot I_{L''},
$$
where $\widetilde {I_{X',L'}}$ is any ideal that
vanishes on $L''$ and restricts on $L'$ to the ideal $I_{X',L'}$ of
$X'$ in $L'$,
and similarly for $\widetilde{I_{X'',L''}}$.
\item{$(b)$} 
$\mu (I_X) = \mu (I_{X',L'})+\mu(I_{X'',L''})+\mu(I_{L'})\mu(I_{L''})$.
\item{$(c)$} 
$\reg(I_X)=\max\set{2, \reg (I_{X'}), \reg (I_{X''})}$.

The simplest way to construct an ideal $\widetilde {I_{X',L'}}$ as
required in part $(a)$ is to choose coordinates in $\PP^r$ so that $x_0,\dots,x_i$
are coordinates for $L'$ while the linear space where $x_0,\dots,x_i$
vanish is a subspace of $L''$. Let $I_{X',L'}$ be the ideal of $X'$ in
the homogeneous coordinate ring of $L'$.  If we write generators for
$I_{X',L'}$ in terms of the coordinates $x_0,\dots,x_i$, then we may
take the same expressions as generators for $\widetilde {I_{X',L'}}$.

Geometrically, this construction amounts to choosing a subspace
$L''_0\subset L''$ complementary to $L=L'\cap L''$, and taking
$\widetilde {I_{X',L'}}$ to be the ideal of the cone with base $X'$
and vertex $L''_0$.  However, this is not the only choice possible: we
may perturb the generators given above by any elements of $I_{L'}\cap
I_{L''}$. This may even change the codimension of $\widetilde
{I_{X',L'}}$. Note that the generators of $I_{X',L'}$ have degree at
least 2. Since $(I_{L'}\cap I_{L''})_{\geq 2} = (I_{L'}\cdot
I_{L''})_{\geq 2}$ such a perturbation will not, as it must not,
change the given form of $I_X$.

\smallskip

{\it Proof of \ref{lin join equations}.\ }  
Since $I_X=I_{X'}\cap I_{X''}$ and $I_{X'}+I_{X''}=I_L$
 there is an exact sequence
$$
0\rTo I_{X}\rTo I_{X'}\oplus I_{X''} \rTo I_{X'\cap X''}\rTo 0.
$$
Since we have assumed that $X\subset \lspan(L'\cup L'')=\PP^r$ is
nondegenerate, the space of linear forms vanishing on $L$ is the 
direct sum of the spaces of linear forms vanishing on $L'$ and $L''$.
If we choose minimal free modules $F', F''$ and surjections
$F'\to I_{L'}$ and $F''\to I_{L''}$ we get a minimal surjection
$F'\oplus F''\to I_L$, and we may write a minimal free resolution
of $I_L$ in the form
$$
\wedge(F'\oplus F'')_{\geq 1}:\qquad
 \cdots\rTo \wedge^2(F'\oplus F'')
\rTo F\oplus F'\rTo I_L\rTo 0.
$$
Let
$\widetilde {I_{X',L'}}$ be any ideal that
vanishes on $L''$ and restricts on $L'$ to the ideal of
$X'$ in $L'$,
and similarly for $\widetilde{I_{X'',L''}}$.
Let $G_1'$ and $G_1''$ be free modules that map 
minimally onto 
$\widetilde {I_{X',L'}}$
and
$\widetilde{I_{X'',L''}}$
respectively, and set $H_1'=F'\oplus G_1'$ and
$H_1''=F''\oplus G_1''$. Because 
$\widetilde {I_{X',L'}}$
reduces modulo $I_{L'}$ to $I_{X',L'}$ the 
induced map $H_1'\to I_{X'}$ is surjective,
and similarly for $H_1''\to I_{X''}$. Choose now 
a minimal free resolution
$$
\HH':\qquad  \cdots \rTo H_2'\rTo H_1'\rTo I_{X'}\rTo 0.
$$
The Koszul complex on the generators of $I_{L'}$
may be written as
$(\wedge F')_{\geq 1}$.
Considering degrees, we see that
each term
$\wedge^i F'$ is a summand of the corresponding term
$H'_i$ of $\HH'$, and we may write 
$H_i'= \wedge^i F'\oplus G_i'$ for a suitable free module $G_i'$.
Define $\HH'', \wedge F''$ and $G_i''$ similarly, and let
$a: H_1'\oplus H_1''\to I_{X'}\oplus I_{X''}$ be the projection.

We will choose a map of complexes
$\phi:\HH'\oplus \HH''\to \wedge(F'\oplus F'')_{\geq 1}$
lifting the natural sum map $I_{X'}\oplus I_{X''}\to I_L$:
$$
\diagram[small,midshaft]
\vdots&&\vdots\cr
\dTo&&\dTo\cr
(G_2'\oplus \wedge^2F')\oplus (G_2''\oplus \wedge^2F'') &\rTo^{\phi_2} &
\wedge^2F'\oplus (F'\otimes F'')\oplus \wedge^2F''\cr
\dTo&&\dTo_\delta\cr
(G_1'\oplus F')\oplus (G_1''\oplus F'') &\rTo^{\phi_1} &F'\oplus F''\cr
\dTo^a&&\dTo_b\cr
I_{X'}\oplus I_{X''} &\rTo &I_{L}.
\enddiagram
$$
We may make this choice so that the restriction of $\phi$
to the subcomplex
$(\wedge F'\oplus \wedge F'')_{\geq 1}\subset \HH$
is the canonical injection into $\wedge(F'\oplus F'')_{\geq 1}$.

The mapping cone of $\phi$ is a complex
$$
\cdots 
\rTo 
\biggl[\wedge^2F'\oplus (F'\otimes F'')\oplus \wedge^2F''\biggr]
\oplus
\biggl[(G_1'\oplus F')\oplus(G_1''\oplus F'')\biggr]
\rTo^{\delta\oplus \phi_1}
F'\oplus F''
$$
whose only homology is $I_X$, occurring at $H_2'\oplus H_2''$. 

Since the map $\phi_1$ is
surjective and $\phi_2$ maps $\wedge^2 F'\oplus \wedge^2F''$ onto
the corresponding summands of 
$\wedge^2(F'\oplus F'')= \wedge^2 F'\oplus (F'\otimes F'')\oplus \wedge^2F''$,
we see that
$I_X$ has a free resolution beginning with
$$
\cdots(G_2'\oplus \wedge^2F')\oplus (G_2''\oplus \wedge^2F'')
\oplus \wedge^3(F'\oplus F'')
\rTo 
F'\otimes F''\oplus
\ker (\phi_1)
\rTo 
I_X
\rTo 
0.
$$
The ideal $I_X$ is embedded in $I_{X'}\oplus I_{X''}$
as the diagonal. Thus 
the image of $F'\otimes F''$ in $I_X$ is computed by
lifting $\delta$ along $\phi_1$ and then composing with the 
map $H'_1\oplus H''_1\to I_{X'}\oplus I_{X''}$, which lands in
the kernel of the projection $I_{X'}\oplus I_{X''}\rTo I_{X'\cap X''}=I_L$. 
Because of our choice of $\phi$, the image of this composite is 
$I_{L'}\cdot I_{L''}$.

Because 
$\widetilde {I_{X',L'}}\subset I_{L''}$ we may
choose $\phi_1{|}_{G'_1}$ to be a map with image
contained in $F''$. We similarly
choose $\phi_1{|}_{G''_1}$ to map into $F'$. It follows
that $\ker(\phi_1)$ is the direct sum of 
$\ker(\phi_1{|}_{G'_1\oplus F''})= {\rm{Transpose}}{(1,-\phi_1{|}_{G_1'})}(G_1')$ and
a corresponding copy of $G_1''$. 
Since $I_X$ is embedded diagonally in $I_{X'}\oplus I_{X''}$
the image of ${\rm{Transpose}}{(1,-\phi_1{|}_{G_1'})}(G_1')$ is
$\widetilde {I_{X',L'}}$. The corresponding result for
$\widetilde {I_{X'',L''}}$ follows symmetrically, proving
the formula $(a)$.

The resolution of $I_X$ that we have just 
constructed is not in general minimal. 
However, the syzygies corresponding to elements of $G_2'\oplus G_2''$
must involve the elements of $G_1'\oplus G_1''$. As these elements
have degrees $\geq 2$, the generators of $G_2'$ and $G_2''$ must
have degrees $\geq 3.$
Thus the map
from $(F'\otimes F'')$ to $S$ sends the generators minimally
onto the generators of $I_{L'}\cdot I_{L''}$. Moreover, the 
generators of $G_1'$ map onto forms that minimally generate the
ideal of $X'$ in its span, and similarly for $X''$. If there 
were a dependence relation of the form 
$$
f'+f''+p=0
$$ 
where
$$
f'\in \Im(G_1'),\ f''\in \Im(G_1'')\ {\rm{and}}\ p\in\Im(F'\otimes F''),
$$
then working modulo the equations of $L'$ we see that 
$f'=0$ and similarly $f''=0$. It follows that $p=0$ as
well, establishing part $(b)$, the
desired relation on minimal numbers of generators.

Statement $(c)$ on the regularity follows at once from the form of 
the above (not necessarily minimal) resolution.
\Box

\corollary {eqns for unions of linears}
Let $X=\bigcup_{i=1}^nL_i\subset \PP^r$ be a
nondegenerate $2$-regular
union of linear subspaces, linearly joined in that order as in 
\ref{main}. If, for each $i\in\{2,\dots,n\}$,
$L'_i$ is a linear complement in $L_i$ for 
$L_{r_i}\cap L_i$, where  
$\lspan(L_1,\dots, L_{i-1})\cap L_{i} = L_{r_i}\cap L_{i}$
with $1\le r_i<i$, then
$$\eqalign{
I_X&= \sum_{j=2}^n 
I_{\lspan(L_j, L_{j+1},\dots,L_{n})}
I_{\lspan(L_1,\dots,L_{j-1}, L_{j+1}',\dots,L_{n}')}\cr
&=
\sum_{j=2}^n 
I_{\lspan(L_j, L_{j+1}',\dots,L_{n}')}
I_{\lspan(L_1,\dots,L_{j-1}, L_{j+1}',\dots,L_{n}').}
}$$

\proof
The case $n=2$ is immediate.  By induction on $n$ we may assume that
the equations of $X'=L_1\cup\dots\cup L_{n-1}$ are given by a
similar formula with $n-1$ in place of $n$. Set $X''=L_n$.  In the
expression for $I_X$ in \ref{lin join equations}, the ideal
$\widetilde{I_{X''/L''}}$ may be taken to be 0.  We may choose
$\widetilde{I_{X',L'}}$ to be the ideal of the cone with base $X'$ and
vertex $L_n'$, as in the remark after \ref{lin join equations}.  Since
taking the ideals of cones commutes with sums and products in an
appropriate sense, and the cone over a span of a collection of linear
spaces is obtained by taking the span with the vertex, we arrive at
the given formula.\Box

\remark{direct sum}
By the above choice of $L_i'$ we have, for each $i<j$, that
$$
\lspan(L_i,L_{i+1}',\dots, L_j')=\lspan(L_i,L_{i+1},\dots, L_j).
$$
Applying this with $i=1, j=n$ we see that $L_1,L_2',\dots, L_{n}'$ 
span $\PP^r$. Moreover, since 
$$
\lspan(L_1,\dots, L_j)\cap L_{j+1} = L_{r_j}\cap L_{j+1}
$$
by hypothesis, we see that $L_{j+1}'$ is disjoint from
$\lspan(L_1,\dots,L_j)$. Thus if we think of $\PP^r$ as lines
in a vector space $V$, then $V$ is the direct sum of spaces
corresponding to $L_1,L_2', \dots, L_{n}'$. It follows that
we may choose variables $x_i$ so
that for each $j$ the space $\lspan(L_{j+1}',\dots,L_{n}')$
is defined by the vanishing of an initial segment
$x_0,\dots, x_{i_j}$, and this set of variables
are coordinates on $\lspan(L_1,\dots,L_j)$. 
If $Y\subset \lspan(L_1,\dots, L_j)$ is a subvariety, then
the cone with base $Y$ and vertex $\lspan(L_{j+1}'\cup \ldots\cup L_{n}')$ is defined
by the equations of $Y$ in $\lspan(L_1,\dots, L_j)$ written in
these coordinates. Moreover, the cone over this variety with
vertex $L_{j+2}'$ is given by the same equations, and so on.

\smallskip

\corollary{products-of-lin-forms} The homogeneous ideal $I\subset S$ of any 
{\rm 2}-regular union of linear spaces is generated by 
products of pairs of distinct (independent) linear forms.\Box

\remark{no inflation} Motivated by the structure in \ref{products-of-lin-forms} 
for the homogeneous ideal of a $2$-regular union of linear spaces one
may believe that such an ideal can always be obtained from the
squarefree monomial ideal of a $2$-regular union of coordinate
subspaces in a larger space by factoring out a sequence of linear
forms that is a regular sequence on every component and any nonempty
mutual intersection of irreducible components. Unfortunately this is
false as the following example shows: Let $X\subset\PP^6$ consist of a
$\PP^3$ and three general $\PP^2$'s sticking out of it, each meeting
the $\PP^3$ in a line, say denoted as $L_i$.  The set $X$ is
small, but the existence a squarefree monomial ideal ``inflation'' as
above would imply that the homogeneous ideal of the union of the three
(general) lines $L_1, L_2$ and $L_3$ in $\PP^3$ would also be
generated by products of pairs of linear forms, which is not the case.

\references
\parindent=0pt 
\frenchspacing 

\item{} M.~Barile, M.~Morales: On the equations defining
minimal varieties, {\sl Comm. Algebra} {\bf 28} (3) (2000), 1223-1239.
\medskip

\item{} M.~Barile, M.~Morales:  On Stanley-Reisner rings of reduction number one, 
{\sl Ann. Scuola Norm. Sup. Pisa Cl. Sci.} (4) {\bf 29} (2000), no. 3, 605--610.
\medskip

\item{} M.~Barile, M.~Morales: On unions of scrolls along linear spaces, 
{\sl Rend. Sem. Mat. Univ. Padova} {\bf 111} (2004), 1--15.
\medskip

\item{} D.~Bayer, D.~Mumford: What can be computed in algebraic geometry?,  
in {\sl "Computational Algebraic Geometry and Commutative Algebra"} 
(D. Eisenbud and L. Robbiano, Eds.), pp. 1--48, Cambridge Univ. Press, 1993.
\medskip

\item{} E.~Bertini: {\sl Introduzione alla geometria proiettiva degli
iperspazi}, Enrico Spoerri, Pisa, 1907.
\medskip

\item{} J.~Blair, B.~Peyton:
An introduction to chordal graphs and clique trees, in 
{\sl Graph theory and sparse matrix computation}, 1--29,
IMA Vol. Math. Appl. {\bf 56}, Springer, New York, 1993. 
\medskip

\item{} G.~Caviglia: Bounds on the Castelnuovo-Mumford regularity
of tensor products, preprint 2003, {\sl Proc. Am. Math. Soc.}, to appear
\medskip

\item{} P.~Del Pezzo: Sulle superficie di ordine $n$ immerse nello spazio di
$n+1$ dimensioni, {\sl Rend. Circ. Mat. Palermo} {\bf 1} (1886).
\medskip

\item{} G.A.~Dirac: On rigid circuit graphs
{\sl Abh. Math. Sem. Univ. Hamburg} {\bf 25} (1-2) (1961), 71--76.
\medskip

\item{} D.~Eisenbud: 
{\sl Commutative Algebra with a View Toward Algebraic Geometry}, 
Springer, New York, 1995. 
\medskip 

\item{} D. Eisenbud: {\sl Geometry of Syzygies}, 
forthcoming book, Springer 2004.
\medskip

\item{} D.~Eisenbud, S.~Goto:
Linear free resolutions and minimal multiplicity. 
{\sl J. Algebra}, {\bf 88} (1984), no. 1, 89--133.
\medskip

\item{} D.~Eisenbud, M.~Green, K.~Hulek, S.~Popescu:
Restricting linear syzygies: algebra and geometry,
preprint 2004.
\medskip

\item{} D.~Eisenbud, J.~Harris:
On varieties of minimal degree (a centennial account), in 
``Algebraic geometry, Bowdoin, 1985 (Brunswick, Maine, 1985)'', 3--13,
{\sl Proc. Sympos. Pure Math.} {\bf 46}, Part 1.
\medskip

\item{} R.~Fr\"oberg: 
Rings with monomial relations having linear resolutions.
{\sl J. Pure Appl. Algebra} {\bf 38} (1985) 235--241.
\medskip

\item{} R.~Fr\"oberg: On Stanley-Reisner rings. In
{\sl ``Topics in algebra, Part 2'' (Warsaw, 1988)}, Banach Center Publ., 
{\bf 26}, Part 2, 57--70, PWN, Warsaw, 1990.
\medskip

\item{} D.R.~Fulkerson, O.A.~Gross: 
Incidence matrices and interval graphs,
{\sl  Pacific J. Math.} {\bf 15} (3) (1965), 835-855.
\medskip

\item{} W.~Fulton: {\sl Intersection Theory}, Springer-Verlag, Berlin, 1984.
\medskip

\item{} W.~Fulton, R.~Lazarsfeld:
Positivity and excess intersection, in 
{\sl ``Enumerative geometry and classical algebraic geometry (Nice, 1981)''}, 
97--105, Progr. Math. {\bf 24}, Birkh\"auser 1982. 
\medskip
 
\item{} R.L.~Graham, P.~Hell: On the history of the minimum spanning tree problem,
{\sl Ann. Hist. Comput.} {\bf 7}, (1985), no. 1, 43--57.
\medskip

\item{} D.~Grayson, M.~Stillman: 
{\sl Macaulay2}: A computer program designed to support computations
in algebraic geometry and computer algebra.  Source and object code
available from {\tt http://www.math.uiuc.edu/Macaulay2/}.
\medskip

\item{} R.~Hartshorne: Complete intersections and connectedness.
{\sl Amer. J. Math.} {\bf 84} (1962), 497--508.
\medskip

\item{} J.~Herzog, T.~Hibi, X.~Zheng: Dirac's theorem on
chordal graphs and Alexander duality, preprint 2003.
\medskip

\item{} R.~Lazarsfeld: {\sl Positivity in Algebraic Geometry I \& II},
 Ergebnisse der Mathematik und ihrer Grenzgebiete, vol {\bf 48} and {\bf 49},
Springer 2004.
\medskip

\item{} Ch.~Miyazaki, W.~Vogel, K.~Yanagawa:
Associated primes and arithmetic degrees,
{\sl J. Algebra} {\bf 192} (1997), no. 1, 166--182.
\medskip

\item{} D.~Mumford: {\sl Lectures on Curves on an Algebraic Surface},
Annals of Math.~Studies {\bf 59},
Princeton University Press, Princeton, New Jersey, 1966.
\medskip

\item{} D.~Mumford: {\sl Algebraic Geometry I: Complex Projective Varieties}, 
Springer, New York, 1976.
\medskip

\item{}
R.C.~Prim: Shortest connection networks and some generalizations, 
{\sl Bell System Tech. J.}, {\bf 36} (1957), 1389-1401.
\medskip

\item{}
J.~Sidman: On the Castelnuovo-Mumford regularity of products of ideal sheaves,
{\sl Adv. Geom.} {\bf 2} (2002), no. 3, 219--229.
\medskip

\item{} R.~Stanley: {\sl Combinatorics and Commutative Algebra}, Second
edition, Progress in Math. {\bf 41}, Birkh\"auser, 1996.
\medskip

\item{} B.~Sturmfels, N.V.~Trung, W.~Vogel: 
Bounds on degrees of projective schemes,  
{\sl Math. Ann.}  {\bf 302} (1995), 417--432.
\medskip

\item{} R.E.~Tarjan:
{\sl Data structures and network algorithms},
CBMS-NSF Regional Conference Series in Applied Mathematics {\bf 44},
Society for Industrial and Applied Mathematics (SIAM), Philadelphia, 1983. 
\medskip

\item{} R.E.~Tarjan, M.~Yannakakis:
Simple linear-time algorithms to test chordality of graphs, test acyclicity of 
hypergraphs, and selectively reduce acyclic hypergraphs,
{\sl  SIAM J. Comput.}  {\bf 13}  (1984),  no. 3, 566--579. 
\medskip

\item{} S.~Xamb\'o:
On projective varieties of minimal degree.
{\sl Collect. Math.} {\bf 32} (1981), no. 2, 149--163. 
\medskip

\vskip .6cm
\widow {.1}
\vbox{\noindent Author Addresses:
\medskip
\noindent{David Eisenbud}\par
\noindent{Department of Mathematics, University of California at Berkeley,
Berkeley CA 94720}\par
\noindent{\tt de@msri.org}
\medskip
\noindent{Mark Green}\par
\noindent{Department of Mathematics, University of California at Los Angeles,
Los Angeles CA 90095-1555}\par
\noindent{\tt mlg@math.ucla.edu}
\medskip
\noindent{Klaus Hulek}\par
\noindent{Institut f\"ur Mathematik, Universit\"at Hannover, D-30060 Hannover,
Germany}\par
\noindent{\tt hulek@math.uni-hannover.de}
\medskip
\noindent{Sorin Popescu}\par
\noindent{Department of Mathematics, SUNY at Stony Brook,
Stony Brook, NY 11794-3651}\par
\noindent{\tt sorin@math.sunysb.edu}\par
}

\end